\documentstyle[12pt,amscd,amsfonts]{amsart}
\newtheorem{thm}[equation]{Theorem}
\numberwithin{equation}{section}
\newtheorem{cor}[equation]{Corollary}

\newtheorem{diag}[equation]{Diagram}

\newtheorem{prop}[equation]{Proposition}

\begin{document}
\raggedbottom
\voffset=-.7truein
\hoffset=0truein
\vsize=8truein
\hsize=6truein
\textheight=8truein
\textwidth=6truein
\baselineskip=18truept
\def\mapright#1{\smash{\mathop{\longrightarrow}\limits^{#1}}}
\def\mapleft#1{\smash{\mathop{\longleftarrow}\limits^{#1}}}
\def\mapup#1{\Big\uparrow\rlap{$\vcenter {\hbox {$#1$}}$}}
\def\mapdown#1{\Big\downarrow\rlap{$\vcenter {\hbox {$\ssize{#1}$}}$}}
\def\mapne#1{\nearrow\rlap{$\vcenter {\hbox {$#1$}}$}}
\def\mapse#1{\searrow\rlap{$\vcenter {\hbox {$\ssize{#1}$}}$}}
\def\mapr#1{\smash{\mathop{\rightarrow}\limits^{#1}}}
\def\lb{[}
\def\ss{\smallskip}
\def\sm{\wedge}
\def\qq{\qquad}
\def\la{\langle}
\def\ra{\rangle}
\def\on{\operatorname}
\def\kbar{{\overline k}}
\def\qed{\quad\rule{8pt}{8pt}\bigskip}
\def\ssize{\scriptstyle}
\def\a{\alpha}
\def\bz{{\bold Z}}
\def\fr{\on{Fr}}
\def\im{\on{im}}
\def\ext{\on{Ext}}
\def\cm{{\cal M}}
\def\sq{\on{Sq}}
\def\eps{\epsilon}
\def\ar#1{\stackrel {#1}{\rightarrow}}
\def\br{\bold R}
\def\bc{\bold C}
\def\si{\sigma}
\def\Ebar{{\overline E}}
\def\Sum{\sum}
\def\tfrac{\textstyle\frac}
\def\tb{\textstyle\binom}
\def\Si{\Sigma}
\def\w{\wedge}
\def\equ{\begin{equation}}
\def\b{\beta}
\def\G{\Gamma}
\def\g{\gamma}
\def\endeq{\end{equation}}
\def\sn{S^{2n+1}}
\def\zp{\bold Z_p}
\def\P{{\cal P}}
\def\cG{{\cal G}}
\def\Hom{\on{Hom}}
\def\ker{\on{ker}}
\def\coker{\on{coker}}
\def\da{\downarrow}
\def\io{\iota}
\def\Om{\Omega}
\def\u{{\cal V}}
\def\e{{\cal E}}
\def\cg{{\cal G}}
\def\exp{\on{exp}}
\def\xbar{{\overline x}}
\def\ebar{{\overline e}}
\def\et{{\widetilde E}}
\def\ni{\noindent}
\def\coef{\on{coef}}
\def\den{\on{den}}
\def\lcm{\on{l.c.m.}}
\def\vi{v_1^{-1}}
\def\ot{\otimes}
\def\psibar{{\overline\psi}}
\def\phibar{{\overline\phi}}
\def\mhat{{\hat m}}
\def\exc{\on{exc}}
\def\ms{\medskip}
\def\ehat{{\hat e}}
\def\dirlim{\on{dirlim}}
\def\tor{\on{Tor}}
\def\psit{\widetilde\psi}
\def\bq{{\bold Q}}
\def\zphat{{\widehat{\bz}_p}}
\def\rhot{{\widetilde\rho}}
\def\vp{v_1^{-1}\pi}
\title[From representation theory to homotopy groups]
{From representation theory to homotopy groups}
\author[D. Davis]{Donald M. Davis}
\address{Lehigh University\\Bethlehem, PA 18015}
\email{dmd1@@lehigh.edu}

\date{}

\maketitle
\section{Introduction}\label{intro}

The $p$-primary $v_1$-periodic homotopy groups of a space $X$, denoted
$\vp_*(X;p)$ or just $\vp_*(X)$,
were defined in \cite{DM}. They are a localization of the actual
homotopy groups, telling roughly the portion which is detected by $K$-theory
and its operations. If $X$ is a compact Lie group or spherically resolved space, each $\vp_i(X;p)$ is a
direct summand of some actual homotopy group of $X$.

By use of a combination of homotopy-theoretic and unstable Novikov spectral
sequence (UNSS)
methods, the groups $\vp_*(X;p)$ were computed 
by the author and coworkers for the following
compact simple Lie groups:
\begin{itemize}
\item $X$ a classical group and $p$ odd (\cite{D});
\item $X$ an exceptional Lie group with $H_*(X;\bold Z)$ $p$-torsion-free
(\cite{BDMi});
\item ($SU(n)$ or $Sp(n)$, $2$) (\cite{BDSU},\cite{BDM});
\item $(G_2,2)$ (\cite{DMG}), ($F_4$ or $E_6$, 3) (\cite{F4}),
and $(E_7,3)$ (\cite{E7}).
\end{itemize}

In \cite{Bous}, Bousfield takes
a new approach to $v_1$-periodic homotopy groups.
He shows that if $X$ is a 1-connected finite $H$-space with
$H_*(X;\bq)$ associative, and $p$ is an odd prime,
then $\vp_*(X;p)$ can be obtained explicitly from
$(K^*(X;\zphat),\psi^p,\psi^r)$, where $r$ is a generator of the group of units
$(\bz/p^2)^\times$. We will review his result in Theorem \ref{Bousthm}.

Let $X$ be any compact simple Lie group, and $p$ an odd prime.
In this paper, we will show how to
compute the second exterior power $\lambda^2$ of generators of
the representation ring $R(X)$, and use this to find a set of generators of
$K^*(X;\zphat)$ on which the Adams operations $\psi^k$ behave in a nice way.
From this, we use Bousfield's theorem to determine $\vp_*(X;p)$.
This approach is totally algebraic. There is no homotopy theory
(except that which went into proving Bousfield's theorem) and no UNSS.

We have used this approach to check the results obtained earlier by
homotopy theory and the UNSS for $(X,p)$ if $X=G_2$, $F_4$, $E_6$, or $E_7$,
and
$p\ge3$, and $X=E_8$ and $p\ge7$. All results are in
agreement, except for one minor mistake in \cite{F4} in $\vp_*(F_4;3)$, which
will be discussed in Section \ref{F4}. Also in Section \ref{F4} we will show
how this new approach resolved two minor matters for $(E_7,3)$ which had been
left unresolved in \cite{E7}.

In this paper, we will focus our
attention on the calculation of $\vp_*(E_8;5)$ and $\vp_*(E_8;3)$, both of which are new. A main impediment
toward finding $\vp_*(E_8;5)$ had been
uncertainty about a product decomposition
which had been claimed by Harper in 1974
in \cite[4.4.1(b)]{Har}. In 1987, Kono questioned Harper's proof, and
Harper agreed to Kono that his proof was flawed. Our methods show that indeed
Harper's claim was incorrect; the asserted product decomposition does not
exist. This will be explained more fully in Proposition \ref{Harprop}.

Our main results are as follows, but we feel that the new methods introduced to
obtain them are of much more interest than the results themselves. Let
$\nu_p(-)$ denote the exponent of $p$ in an integer.
\begin{thm} \label{main}Let $\la m,k\ra_r=\min(\nu_5(m-k)+r,k)$. Then
\begin{eqnarray*}\vp_{2m}(E_8;5)&\approx& \vp_{2m-1}(E_8;5)\\
&\approx&
\begin{cases}0&\text{if $m$ is even}\\
\bz/5^{\max(\la m,7\ra_4,\la m,11\ra_4,\la m,19\ra_4,\la m,23\ra_4)}&
\text{if $m\equiv 3$ mod $4$}.\end{cases}\end{eqnarray*}
If $m\equiv1$ mod $4$, then
$\vp_{2m}(E_8;5)\approx \bz/5^e\oplus\bz/5$, and $\vp_{2m-1}(E_8;5)\approx
\bz/5^{e+1}$, where
$$e=\max(\la m,13\ra_2,\la m,17\ra_2,\la m,29\ra_2).$$
\end{thm}

\begin{thm}\label{e83thm} If $i\equiv3,4$ mod $4$, then $\vp_i(E_8;3)=0$.
For any integer $k$, 
$$\vp_{4k+1}(E_8;3)\approx \vp_{4k+2}(E_8;3)\approx 3^e,$$
where
$$e=\begin{cases}\min(7+\nu_3(k-9-3^{13}),21)&\text{if $k\equiv0$ mod $9$}\\
6&\text{if $k\equiv1$ mod $3$}\\
\min(9+\nu_3(k-11),24)&\text{if $k\equiv2$ mod $9$}\\
\min(7+\nu_3(k-6-2\cdot3^7),15)&\text{if $k\equiv3,6$ mod $9$}\\
\min(10+\nu_3(k-14),30)&\text{if $k\equiv5$ mod $9$}\\
\min(9+\nu_3(k-8),18)&\text{if $k\equiv8$ mod $9$.}
\end{cases}$$
\end{thm}

This completes the computation of $\vp_*(X;p)$ for all compact simple Lie groups $X$ and odd primes $p$, a project which was suggested to the author by Mimura in 1989. The situation when $p=2$, which was part of Mimura's suggested project, is much more delicate. The author hopes to be able to adapt Bousfield's theorem to the prime 2, but that work is still in very preliminary stages.

There are many computations in this project which would be intractable to do by
hand. Specialized software
{\tt LiE} (\cite{Lie}) is used to determine the second
exterior power operations in $R(E_8)$. A nontrivial algorithm was required to
get this information into the form of an $8\times 8$ matrix
of integers, some of them 16 digits long, which can be interpreted as giving
$\psi^2$ on a canonical basis of 
$PK^1(E_8)$, the primitive elements. This portion of the work will
be described in Section \ref{rep}.

The eigenvalues of this matrix
are $2^e$ for $e\in R=\{1,7,11,13,17,19,23,29\}$,
corresponding to the fact that rationally $E_8$ is equivalent to
$\prod_{e\in R}S^{2e+1}$. Using {\tt Maple}, we find the associated
eigenvectors.
The determinant  of the matrix of these eigenvectors is
$$D=2^{61}3^{32}5^{10}7^911^413^417^319^223^229.$$
This implies that localized at a prime $p$
greater than 29, these eigenvectors span
$PK^1(E_8)_{(p)}$, which then is isomorphic to $PK^1(\prod_RS^{2e+1})_{(p)}$ as
a
module over all Adams operations $\psi^k$, and hence by Bousfield's theorem has
$$\vp_*(E_8;p)\approx\vp_*(\prod_R S^{2e+1};p).$$

If, for example, $p=29$, we can find two of the eight eigenvectors, $v$ and $w$, for which
$v':=(v-w)/29$ is integral. The set of vectors obtained from the eight
eigenvectors by replacing $v$ by $v'$ has its determinant equal to $D/29$,
which is a unit in $\bz_{(29)}$, and so this set spans $PK^1(E_8)_{(29)}$.
The eigenvectors $v$ and $w$ correspond to eigenvalues $2^1$ and $2^{29}$.
Then $\psi^k(v)=kv$ and $\psi^k(w)=k^{29}w$ for all integers $k$, 
and so we can determine $\psi^k(v')$, and
from this use Bousfield's theorem to find $\vp_*(E_8;29)$, which agrees with that deduced in \cite{BDMi} from the decomposition
$$E_8\simeq_{29} B(3,59)\times S^{15}\times S^{23}\times S^{27}\times
S^{35} \times S^{39}\times S^{47}.$$
The point is that the value of the determinant $D$, computed blindly from
representation theory, is intimately related to the 
decomposition of $E_8$ when localized at each prime. Note also that the matrix analysis
shows that the portion of $K^1(E_8;\bq)$ which must be modified to pass
to $K^1(E_8;\bz_{(29)})$ is the portion related to $S^3$ and $S^{59}$,
consistent with the homotopy analysis.

Because of the $5^{10}$ factor in $D$, we must 10 times replace vectors by $\frac15$ times a difference of vectors
in order to find a set of vectors whose determinant is a unit
in $\bz_{(5)}$.
On this set, an explicit formula for the Adams operations $\psi^k$ can be
given. This portion of the work will be described in Section \ref{vecs}.
We also show there how performing these basis changes for all relevant primes enables a 
determination of the Adams operations $\psi^k$ in $K^*(X)$ 
(not localized at a prime) for all exceptional Lie groups $X$ and all $k$.

The Adams operation formulas are of the sort that allow us to draw inferences about attaching maps in the localized Lie groups. This new homotopy-theoretic information has been derived here just from our
representation-based calculations together with Adams' $e$-invariant work
(\cite{JXIV}). This will be discussed in Sections \ref{odd} and \ref{3calc}.

We feed this information into Bousfield's theorem, which, after a good
deal of manipulation, yields the results for $\vp_*(E_8;5)$.
The computations for $\vp_{2m}(E_8)$ are given in Section \ref{even}, and those
for $\vp_{2m-1}(E_8)$ are given in Section \ref{odd}, which also includes some
useful general results, such as periodicity of the number of summands, and the use of exact sequences.
In Section \ref{3calc}, we perform a similar analysis to obtain the result in Theorem \ref{e83thm} for $\vp_*(E_8;3)$.

Bendersky and Thompson (\cite{BT}) have recently constructed an unstable
Bousfield-Kan spectral sequence based on $K_*K$, where $K$ represents
periodic $K$-theory. It possesses some advantages over the 
$BP$-based UNSS used in papers such as \cite{BDMi} and \cite{E7},
especially in that $K_*(E_8)_{(p)}$ is a free
commutative algebra, whereas $BP_*(E_8)$ is not. Using the result obtained in our Proposition \ref{niceops}
about the Adams operations in $K^*(E_8)_{(5)}$,
which effectively implies that there is an
$\alpha_3$ attaching map between cells which would have been in
separate factors if Harper's asserted splitting had been correct, together
with their $K_*K$-based spectral sequence, Bendersky and Thompson are
able to compute
$\vp_*(E_8;5)$ in a manner which is arguably more insightful than the
computation here. Their
method requires information about the homotopy theory of
$E_8$ (e.g. Steenrod operations),
while ours requires information about the representation theory.
But so far homotopy theory has been unable to provide the complete picture
(i.e., the $\alpha_3$ attaching map), for which it had to rely on the
representation-theoretic approach presented here.
Moreover, our result for $\vp_*(E_8;3)$ still seems totally inaccessible to UNSS-type methods.

The author would like to thank Martin Bendersky, Pete Bousfield, and
Mamoru Mimura for helpful comments on this project.

\section{Representation theory and $\psi^2$ in $K$-theory}\label{rep}
In this section, we use representation theory to determine the Adams operation
$\psi^2$ in $K^*(E_8)$, and present an algorithm by which this can be done for
any compact simple Lie group.

Let $G$ be a simply-connected compact Lie group. Bousfield's approach to
$\vp_*(G;p)$ (\cite{Bous})
requires as input certain Adams operations on the primitives
$PK^1(G;\zphat)$. Bousfield suggested to the author the relationship with
exterior powers in the representation ring $R(G)$ described in the next two
paragraphs.

Let $I$ denote the augmentation ideal in $R(G)$. Hodgkin's theorem (\cite
{Hod}) implies that there is an isomorphism
\begin{equation}\label{II} I/I^2\to PK^{-1}(G),\end{equation}
which may be viewed either as induced from the composition $R(G)\to
K^0(BG,*)\to K^{-1}(G)$ or from Hodgkin's function $\b$ that views a
representation $\rho:G\to U(n)$ as a homotopy class in $[G,U]=K^{-1}(G)$.
Although Hodgkin doesn't write the isomorphism (\ref{II}), he describes $R(G)$
in such a way that $I/I^2$ is clearly the free abelian group on the reduced basic
representations $\rhot_1,\ldots,\rhot_l$, and shows that $PK^{-1}(G)$ is the free
abelian group on $\b(\rhot_1),\ldots,\b(\rhot_l)$. The simple description of
$\b$ makes it clear that (\ref{II}) respects the exterior power operations
$\lambda^n$.

Adams operations are related to exterior powers by the Newton formula
$$\psi^n(a)-\lambda^1(a)\psi^{n-1}(a)
+\cdots+(-1)^{n-1}\lambda^{n-1}(a)\psi^1(a)+
(-1)^nn\lambda^n(a)=0,$$
which implies that $\psi^n=(-1)^{n+1}n\lambda^n$ in $I/I^2$. By \cite[5.3]{VF},
$\psi^n$ in $K^1(G)$ corresponds to $\psi^n/n$ in $K^{-1}(G)$. Thus
$\psi^n$
in $PK^1(G)$ corresponds to $(-1)^{n+1}\lambda^n$ in $I/I^2$.

Now let $G$ be a compact simple Lie group of rank $l$ (e.g., $E_8$ of rank 8).
The representation theory of $G$ is equivalent to that of the associated Lie
algebra $\cg$. Associated to $\cg$ is a set $\Lambda$ of weights, a subset
$\Lambda^+$ of dominant weights, and a subset $\{\lambda_1,\ldots,\lambda_l\}
\subset\Lambda^+$ such that $\Lambda$ (resp. $\Lambda^+$) is the free abelian
group (resp. free abelian monoid) generated by $\lambda_1,\ldots,\lambda_l$. 
(e.g., \cite[p.67]{Hum}.) The set $\Lambda$ is given a partial order by
$\sum m_i\lambda_i\le\sum m_i'\lambda_i$ if and only if $m_i\le m_i'$ for
$1\le i\le l$.

To each irreducible representation of $\cg$ is associated a finite set of
weights with multiplicities.
 It is a theorem that
one of these weights is larger than all the others, and it occurs with
multiplicity 1. (e.g., \cite[\S 21.1]{Hum}.) This highest weight is
dominant. The \lq\lq highest weight'' defines a function from the set of
isomorphism classes of irreducible representations of $\cg$ to $\Lambda^+$, and this function is bijective. It is a theorem
(See, e.g., \cite[3.3]{Hod}) that $R(\cg)$ is a polynomial algebra generated by
the irreducible representations $\rho_1,\ldots,\rho_l$ which have
$\lambda_1,\ldots,\lambda_l$ as highest weights.

If ${\bold m}=(m_1,\ldots,m_l)$ is an $l$-tuple of nonnegative integers, let
$V(\bold m)$ denote the unique irreducible representation with highest weight
$m_1\lambda_1+\cdots+m_l\lambda_l$.
We will need three types of information about representations.
\begin{itemize}
\item The dimension of $V(\bold m)$ (as a complex vector space).
\item The second exterior power $\lambda^2(V(\bold m))$, expressed as
$\sum c_jV({\bold k}_j)$ for nonnegative integers $c_j$ and
$l$-tuples ${\bold k}_j$ of nonnegative integers.
\item The tensor product $V({\bold m})\ot V(\bold n)$, expressed as $\sum c'_j
V({\bold k}'_j)$.
\end{itemize}

There are algorithms for each of these, implemented conveniently in the
software {\tt LiE} (\cite{Lie}). For $\dim(V(\bold m))$, Weyl's formula
(\cite[p.139]{Hum}) is used. For $\lambda^2(V(\bold m))$, a formula of
\cite{And} for symmetrized products is used. For $V({\bold m})\ot V(\bold n)$,
Klimyk's formula (\cite{Kli}) is used.

Let $\rhot_i=\rho_i-\dim(\rho_i)$.  We need $\lambda^2(\rhot_i)$ as
a linear combination of $\rhot_1,\ldots,\rhot_l$ in $I/I^2$.
The software gives us $\lambda^2(\rho_i)$ as $\sum c_jV({\bold k}_j)$, and hence
\begin{eqnarray}\lambda^2(\rhot_i)
&=&\lambda^2(\rho_i)+\lambda^1(\rho_i)\lambda^1(-\dim(\rho_i))
+\lambda^2(-\dim(\rho_i))\nonumber\\
&=&\sum c_jV({\bold k}_j)-\dim(\rho_i)\rho_i+
\tbinom{-\dim(\rho_i)}2.\label{dimfor}\end{eqnarray}
We can iterate
tensor product computations to write monomials $\rho_1^{e_1}\cdots\rho_l^{e_l}$
as linear combinations of $V(\bold k)$'s, and then apply an easy row reduction
to this result to
write $V(\bold k)$'s as polynomials in the $\rho_i$'s, or, after
manipulating polynomials,
in the $\rhot_i$'s. In $I/I^2$, we ignore products of $\rhot_i$'s. Substituting
the formulas for $V({\bold k}_j)$ as linear combinations of $\rhot_i$'s into
(\ref{dimfor}) yields the desired expression of
$\lambda^2(\rhot_i)$ as a linear combination of $\rhot_j$'s in $I/I^2$.

We illustrate how this works in the simple example of the exceptional Lie group
$G_2$, and then show how the computations can be expedited. The software tells
us $\lambda^2(\rho_1)=\rho_2+\rho_1$ and $\lambda^2(\rho_2)=V(3,0)+\rho_2$.
(Remember, $\rho_2=V(0,1)$.) Also, $\dim(\rho_1)=7$ and $\dim(\rho_2)=14$.
Thus
$$\lambda^2(\rhot_1)=\rho_1+\rho_2-7\rho_1+\tbinom{-7}2=-
6(\rhot_1+7)+(\rhot_2+14)+\tfrac{7\cdot8}2=-6\rhot_1+\rhot_2.$$
To find $\lambda^2(\rhot_2)$, we need to express $V(3,0)$ as a polynomial in
$\rho_1$ and $\rho_2$. The software tells us
\begin{eqnarray*}\rho_1\ot \rho_1&=&V(2,0)+\rho_1+\rho_2+1\\
\rho_1\ot\rho_2&=&V(1,1)+V(2,0)+\rho_1\\
\rho_1\ot V(2,0)&=&V(3,0)+V(1,1)+V(2,0)+\rho_2+\rho_1.
\end{eqnarray*}
This allows us to compute the second equation in
$$\lambda^2(\rho_2)=V(3,0)+\rho_2=\rho_1^3-2\rho_1\rho_2-\rho_1^2-\rho_1,$$
from which we derive
\begin{eqnarray}\lambda^2(\rhot_2)&=&
(\rho_1^3-2\rho_1\rho_2-\rho_1^2-\rho_1)-14\rho_2+\tbinom{-14}2\nonumber\\
&=&{\rhot_1}^3+20{\rhot_1}^2+104\rhot_1-2\rhot_1\rhot_2-28\rhot_2\nonumber\\
&\equiv&104\rhot_1-28\rhot_2\text{ mod }I^2.\label{I2}\end{eqnarray}
By our earlier remarks, this implies that $PK^1(G_2)$ has basis $\{x_1,x_2\}$
with $\psi^2(x_1)=6x_1-x_2$ and $\psi^2(x_2)=-104x_1+28x_2$.

This procedure can be expedited by just looking at linear terms.
If $f$ is an element of $R(G)$, let $L(f)$ denote the first-order terms when
$f$ is written as a polynomial in $\rhot_1,\ldots,\rhot_l$.
We have
\begin{equation}\label{14}L(\lambda^2(\rhot_2))=L(V(3,0)+\rho_2)-
14\rhot_2,\end{equation}
and from the tensor product equations above, a type of differentiation yields
\begin{eqnarray*}
2\dim(\rho_1)\rhot_1&=&L(V(2,0))+\rhot_1+\rhot_2\\
\dim(\rho_2)\rhot_1+\dim(\rho_1)\rhot_2&=&L(V(1,1))+L(V(2,0))+\rhot_1\\
\dim(V(2,0))\rhot_1+\dim(\rho_1)L(V(2,0))&=&L(V(3,0))+L(V(1,1))\\
&&+L(V(2,0))+\rhot_2+\rhot_1.
\end{eqnarray*}
The first equation says $L(V(2,0))=13\rhot_1-\rhot_2$, then the second says
$L(V(1,1))=8\rhot_2$, and then, using $\dim(V(2,0))=27$, the third equation
says $L(V(3,0))=104\rhot_1-15\rhot_2$, which when substituted into (\ref{14})
yields (\ref{I2}).

The {\tt LiE} program that
implements this expedited algorithm for $E_8$ is listed
and described in Section \ref{pgm}. The one subtlety is how to know which
tensor products to compute. The dominant weights are ordered by
height\footnote{sometimes called \lq\lq level''}, which is the sum of the
coefficients when they are written as roots. For example, in $\cg_2$ the
weights $\rho_1$ and $\rho_2$ correspond to roots $3\a_1+2\a_2$ and
$2\a_1+\a_2$, respectively, and so the height of $V(m_1,m_2)$ is $5m_1+3m_2$.
An elementary result states that $V({\bold m})\ot V({\bold n})=V(\bold m+\bold
n)+$
terms of height less than that of $V(\bold m+\bold n)$. For every term
$V(\bold l)$ with $\sum l_i>1$ which occurs in $\lambda^2(\rho_i)$, we
choose a way of writing $\bold l=\bold m+\bold n$, and differentiate the
equation
$$V({\bold m})\ot V({\bold n})=V({\bold l})+\text{terms of lower height}$$
to inductively obtain formulas for $L(V(\bold l))$. It can happen that
$V({\bold m})\ot V(\bold n)$
might contain terms $V(\bold k)$ which did not appear
in $\lambda^2(\rho_i)$. If so, we also find $L(V(\bold k))$ by the same
method, i.e., differentiating a formula for $V({\bold a})\ot V(\bold b)$ where
$\bold a+\bold b=\bold k$.

When this algorithm is performed for $E_8$, we obtain the matrix (\ref{rhomat})
for $\lambda^2$ on the basis $\{\rhot_1,\ldots,\rhot_8\}$ of $I/I^2$. Thus,
for example,
$$L(\lambda^2\rhot_1)=-3628\rhot_1-\rhot_2+\rhot_3+3875\rhot_8.$$

\begin{scriptsize}
\begin{eqnarray}
\left(\begin{matrix}
-3628&1829621&12625838007&-1270362010619556\\
-1&-116621&-146298269&18170270443687\\
1&-496&-5835130&582917207249\\
0&1&3875&-468700376\\
0&0&-177629&26815340999\\
0&-247&12587859&-2027479372896\\
0&150877&-392633383&78837408033778\\
3875&-4549375&10807381790&-2860474034106800\\
\end{matrix}
\right.\nonumber\\
\left.
\begin{matrix}
2706011993074&-401581533&0&0\\
-40039592220&6661497&0&0\\
-1242615998&185628&0&0\\
1073250&-249&0&0\\
-68699627&27001&-1&0\\
5393300762&-2538745&248&0\\
-233942373952&156457497&-30876&1\\
9317251205935&-7999393170&2573495&-247
\end{matrix}
\right)\label{rhomat}
\end{eqnarray}
\end{scriptsize}

We obtain the following important corollary of this computation, where
$\b$ is the isomorphism of (\ref{II}).
\begin{cor}\label{psi2} With respect to the basis $\{\b(\rhot_1),\ldots,\b(\rhot_8)\}$ of $PK^1(E_8)$,
the matrix of $-\psi^2$ is given by $($\ref{rhomat}$)$.
\end{cor}

\section{Nice form for $\psi^k$ in $PK^1(E_8)_{(5)}$ and $PK^1(X)$}\label{vecs}
If $p=3$ or 5, then 2 generates $(\bz/p^2)^\times$, and so Bousfield's theorem
requires knowledge of $\psi^2$ and $\psi^p$. A computation of $\psi^p$ similar
to that of the previous section could be made (provided enough computer time
and space is available), and the results (Corollary \ref{psi2} and its
analogue) plugged into Bousfield's theorem to give results for $\vp_*(E_8;5)$.
However, the matrix (\ref{rhomat}) is so unwieldy that
it would be very difficult to obtain
a nice form for the resulting groups. For this reason, we find
a new basis of
$PK^1(E_8)_{(5)}$ on which the action of $\psi^2$ has a nicer form.
That is the purpose of this section. Moreover, as we shall see, the new basis
will be one on which every $\psi^k$ can be determined
at the same time as $\psi^2$. We will also perform similar computations
for Adams operations in $K(X)$ (unlocalized) for all exceptional Lie
groups $X$.

To this end, we use {\tt Maple} to find the eigenvalues and eigenvectors of
a matrix $M$, which is defined to be the negative of the matrix
(\ref{rhomat}). This is the matrix of $\psi^2$ on $PK^1(E_8)$.
We are not surprised to find that the eigenvalues of $M$
are $2^1$, $2^7$, $2^{11}$, $2^{13}$, $2^{17}$, $2^{19}$, $2^{23}$, and
$2^{29}$ because of the rational equivalence
\begin{equation}\label{spheres}
E_8\simeq_{\bq}S^3\times S^{15}\times S^{23}\times S^{27}\times S^{35}
\times S^{39}\times S^{47}\times S^{59}\end{equation}
and the fact that $\psi^k$ acts as multiplication by $k^n$ on $K^1(S^{2n+1})$.

A matrix whose columns are eigenvectors
of $M$ corresponding to the eigenvalues listed above in increasing order is

\begin{scriptsize}
\begin{eqnarray}
\left(\begin{matrix}
418105625&451155607289497&-3133156733386433&-2595116726135\\
4168750&3797965233710&16166554278770&49280463350\\
23125&20720212181&97338188051&212788925\\
1&873857&5050967&7697\\
377&326320702&2019676642&2670694\\
119249&99339310201&777470688031&532675273\\
27753998&16454843873197&398907853660267&-470842549139\\
9022308750&-26386060414578330&-79539474507550230&68079530157990
\end{matrix}
\right.\nonumber\\
\left.
\begin{matrix}
-7044348025&-114112039&17026841&2691065\\
-3895749830&-1864130&-136850&-38570\\
21836035&56893&-7787&-1235\\
1471&1&1&1\\
635402&206&446&-58\\
294773639&-86407&-48007&4409\\
-42760408957&12566861&3586541&-174307\\
3836612952570&-810585690&-191313690&6425670
\end{matrix}
\right)\label{evecs}
\end{eqnarray}
\end{scriptsize}

\noindent
The numbers in the columns are coefficients with respect to the basis
of Corollary \ref{psi2}.

The determinant of this matrix is a 68-digit integer which factors as
\begin{equation}\label{bigg}2^{61}3^{32}5^{10}7^911^413^417^319^223^229
\end{equation}
This determinant gives a lot of information. First, it says that localized at a
prime greater than 29, the eigenvectors form a basis, consistent with the known
result that $E_8$ localized at such primes is equivalent to a product of
spheres. Let $$\{v_1,v_7,v_{11},v_{13},v_{17},v_{19},v_{23},v_{29}\}$$
denote the columns of (\ref{evecs}).
It is very important to note that since these vectors satisfy $\psi^2(v_i)=2^i
v_i$, they correspond to the sphere factors in (\ref{spheres}), and hence also
satisfy
\begin{equation}\psi^k(v_i)=k^iv_i\label{psisph}\end{equation}
for any positive integer $k$.

For primes $p$ satisfying $11\le p\le 29$, it was shown in \cite{MNT}
that $E_8$ is $p$-equivalent to a certain product of spheres and sphere bundles
over spheres with $\a_1$ attaching maps. For such primes, the number of sphere
bundles equals the exponent of $p$ in (\ref{bigg}). We use $p=23$ to illustrate
how combinations of the eigenvectors in (\ref{evecs}) correspond to these
product decompositions, providing somewhat more detail than we did in our brief sketch for $p=29$ in Section \ref{intro}.

 We note that $v_1':=(v_1-v_{23})/23$ and $v_7':=
(v_7-18v_{29})/23$ are integral. If $v_1$ and $v_7$ are replaced by $v_1'$ and
$v_7'$ in (\ref{evecs}), then the determinant is divided by $23^2$, and so the
new set of vectors is a basis for $PK^1(E_8)_{(23)}$. It follows from
(\ref{psisph})
that $\psi^k(v_1')=kv_1'+\frac1{23}(k-k^{23})v_{23}$ and
$\psi^k(v_7')=k^7v_7'+\frac{18}{23}(k^7-k^{29})v_{29}$. This agrees with the determination of $\psi^k$
in sphere bundles $B(3,47)$ and $B(15,59)$ with attaching maps $a_1$ given in \cite{JXIV}. Thus as
Adams modules
$$PK^1(E_8)_{(23)}\approx PK^1(B(3,47)\times B(15,59)\times S^{23}\times S^{27}
\times S^{35}\times S^{39})_{(23)}.$$
By Bousfield's theorem, these two spaces will have isomorphic $v_1$-periodic
homotopy groups. Of course, we already knew that by the product
decomposition of \cite{MNT},
but here we are getting it without relying on the \cite{MNT} result.
In \cite{BDMi}, $v_1$-periodic homotopy groups of $\a_1$ sphere bundles over
spheres were determined by the UNSS, and $\vp_*(E_8;p)$ deduced for $p\ge11$
using the product decomposition of \cite{MNT}. In Proposition \ref{sphbdl}, we will
determine the $v_1$-periodic homotopy groups of these sphere bundles by Bousfield's
theorem, giving us a self-contained computation. This can be done
for $(E_8,p)$ for any prime $p\ge11$.

When $p=5$,
we use {\tt Maple} to help us find 10 combinations of vectors that are
divisible by 5. These are
\begin{eqnarray*}v_1'&=&(v_1+2v_{13})/5\\
v_{13}'&=&(v_{13}+3v_{17})/5\\
v_{17}'&=&(v_{17}-v_{29})/5\\
v_{13}''&=&(v_{13}'+v_{17}'-v_{29})/5\\
v_7'&=&(v_7-v_{11})/5\\
v_{11}'&=&(v_{11}-2v_{19})/5\\
v_{19}'&=&(v_{19}-v_{23})/5\\
v_7''&=&(v_7'-v_{11}'+2v_{19}')/5\\
v_{11}''&=&(v_{11}'+v_{19}'+2v_{23})/5\\
v_7'''&=&(v_7''-v_{19}'-2v_{23})/5\end{eqnarray*}
The way that these turn out to be grouped, with 1, 13, 17, and 29 related in
one group, and 7, 11, 19, and 23 related in the other group, is consistent with
Wilkerson's product decomposition (\cite{Wil}) of ${E_8}_{(5)}$ as a product of
two spaces whose rational types correspond to these two groupings.
The matrix $(v_1',v_{13}'',v_{17}',v_{29},v_7''',v_{11}'',v_{19}',v_{23})$
has determinant a unit in $\bz_{(5)}$, and so its columns form a basis
for $PK^1(E_8)_{(5)}$. We rename the classes $(x_1,x_{13},x_{17},x_{29},x_7,
x_{11},x_{19},x_{23})$ and compute $\psi^k(x_i)$ using (\ref{psisph}).
 We obtain the following result.
\begin{prop}\label{niceops} $PK^1(E_8)_{(5)}$ has basis
 $\{x_1,x_{13},x_{17},x_{29},x_7,x_{11},x_{19},x_{23}\}$ satisfying
\begin{eqnarray*}
\psi^k(x_1)&=&k^1x_1+10k(k^{12}-1)x_{13}-8k(k^{12}-1)x_{17}+\tfrac45k(k^{12}-
1)x_{29}\\
\psi^k(x_{13})&=&k^{13}x_{13}+\tfrac45k^{13}(k^4-1)x_{17}+\tfrac2{25}k^{13}
(1+2k^4-3k^{16})x_{29}\\
\psi^k(x_{17})&=&k^{17}x_{17}-\tfrac15k^{17}(k^{12}-1)x_{29}\\
\psi^k(x_{29})&=&k^{29}x_{29}\\
\psi^k(x_7)&=&k^7x_7+\tfrac25k^7(1-k^4)x_{11}+\tfrac1{25}k^7(3-2k^4-
k^{12})x_{19}\\
&&+\tfrac1{125}k^7(32+16k^4-k^{12}-47k^{16})x_{23}\\
\psi^k(x_{11})&=&k^{11}x_{11}+\tfrac15k^{11}(1-k^8)x_{19}+\tfrac1{25}k^{11}
(9k^{12}-k^8-8)x_{23}\\
\psi^k(x_{19})&=&k^{19}x_{19}+\tfrac15k^{19}(1-k^4)x_{23}\\
\psi^k(x_{23})&=&k^{23}x_{23}\end{eqnarray*}
\end{prop}
A formula $\psi^k(x_n)=k^nx_n+\tfrac u5k^n(k^{4m}-1)x_{n+4m}+\cdots$, with
 $u$ a unit in $\bz_{(5)}$ and $m\not\equiv0$ mod 5, is of the type that would be obtained if
the space had cells of dimension $2n+1$ and $2n+1+8m$ with attaching map
$\a_m$. The contributions of these cells to $\vp_*(X)$ will be as if they
had the $\a_m$ attaching map. The formula for $\psi^k(x_{17})$ is, at the very
least, strongly suggestive that in ${E_8}_{(5)}$ the 35- and 59-cells are
connected by $\a_3$. This would contradict the product decomposition asserted
in \cite[4.4.1b]{Har}, which said that one Wilkerson factor could be further
decomposed as $X(3,59)\times X(27,35)$. One could make the argument more
precise either by noting that it is impossible that a space whose Adams
operations can be written as in Proposition \ref{niceops} can be decomposed in
this way
(i.e., no change of basis can split the Adams operations), or by noting
that $\vp_*(E_8;5)$, as computed in the next section, is incompatible with
such a decomposition. Thus we obtain the following result.
\begin{prop}\label{Harprop} In Wilkerson's decomposition $($\cite[2.3]{Wil}$)$
of $(E_8)_{(5)}$ as $X_0\times X_2$, both factors are indecomposable.
In particular, the product decomposition of $X_2$ asserted in \cite[4.4.1b]{Har} is not valid.
\end{prop}

By performing, for all relevant primes, changes of basis of the sort illustrated above for $E_8$ when $p=5$ or 23, we obtain, for each exceptional Lie group $X$, bases for $PK^1(X)$ on which we can compute
$\psi^k$. We obtain the following results. In all of them, let $B_i=\b(\rhot_i)$, where $\b$ and $\rhot_i$ are as in Section \ref{rep}.

\begin{prop}\label{psiG2}
A basis for $PK^1(G_2)$ is given by $\{y_1,y_5\}$ with $y_1=B_1$ and
$y_5=-4B_1+B_2$. For all integers $k$,
\begin{eqnarray*}
\psi^k(y_1)&=&ky_1+\tfrac1{30}(k-k^5)y_5\\
\psi^k(y_5)&=&k^5y_5\end{eqnarray*}
\end{prop}

The nice feature of this result and the subsequent ones for the other exceptional Lie groups is that they are a result about integral 
$K$-theory (i.e., not localized at a prime) and the Adams operations
have a nice form (triangular, among other features). Since $PK^1(X)$
generates $K^*(X)$ for compact simple Lie groups $X$, multiplicativity of the Adams operations allows us to deduce the Adams operations on all
of $K^*(X)$. Note that the classes $y$ are subscripted by the exponent
$e$ such that $\psi^k(y)=k^ey+$ other terms.

The method of proof in each case is to 
\begin{itemize}
\item use {\tt LiE} to obtain the matrix of $\psi^2$ on
$\{B_1,\ldots,B_l\}$ similarly to (\ref{rhomat}),\\
\item use {\tt Maple} to find
eigenvectors of this matrix similarly to (\ref{evecs}), and
note that these vectors satisfy $\psi^k(v_i)=k^iv_i$ for all $k$,\\  \item use {\tt Maple} to repeatedly replace
vectors $v$ by $(v-w)/p$, where $p$ is a prime which divides the determinant of the matrix of vectors and $w$ is a linear combination
of vectors which appear after $v$ in the most recent set of vectors, until the determinant is $\pm1$, and\\
\item use {\tt Maple} and (\ref{psisph})
to compute $\psi^k$ on the final set of vectors, since they are
explicit combinations of the eigenvectors.
\end{itemize}

In \cite{W1} and \cite{W2}, Watanabe computed the Chern character on a
certain set of generators of $K(G_2)$, $K(F_4)$, and $K(E_6)$, and in 
\cite{W3} he explained the well-known way in which this would allow
one to determine the Adams operations. Using {\tt LiE}, his set of generators 
can be expressed as linear combinations of ours, and the results for $\psi^k$ 
which could be read off from his results for $ch$
can be related to ours. We checked this for $G_2$ and $F_4$ and found the results to be 
in agreement. However, it should be pointed out that
his generators do not have the nice feature of having a triangular matrix for $\psi^k$. 
Also, in June 1998, Watanabe told the author that he felt that his methods would not work 
for $E_7$ and $E_8$. An internal check on our results, which was performed for each
exceptional Lie group, is to transform the $\psi^k$ formulas to the
basis of $B_i$'s, and then set $k=2$ and compare with the formulas given by {\tt LiE}.

Now we state the results for the other exceptional Lie groups.
\begin{prop}\label{psiF4} A basis for $PK^1(F_4)$ is given by
$\{y_1,y_5,y_7,y_{11}\}$ with $y_1=-B_4$, $y_5=B_1-3B_4$, $y_7=-2B_1+B_3-15B_4$, and $y_{11}=
-6B_1+B_2-11B_3+102B_4$. For all integers $k$,
\begin{eqnarray*}\psi^k(y_1)&=&ky_11\tfrac1{10}(k-k^5)y_5+(-\tfrac1{70}k+\tfrac1{120}k^5+
\tfrac1{168}k^7)y_7\\
&&+(-\tfrac1{4620}k+\tfrac1{6720}k^5+\tfrac1{13440}k^7-\tfrac1{147840}k^{11})y_{11}\\
\psi^k(y_5)&=&k^5y_5+\tfrac1{12}(k^5-k^7)y_7+(\tfrac1{672}k^5-\tfrac1{960}k^7-\tfrac1{2240}k^{11})y_{11}\\
\psi^k(y_7)&=&k^7y_7+\tfrac1{80}(k^7-k^{11})y_{11}\\
\psi^k(y_{11})&=&k^{11}y_{11}
\end{eqnarray*}
\end{prop}

\begin{prop}\label{psiE6} A basis for $PK^1(E_6)$ is given by 
$\{y_1,y_4,y_5,y_7,y_8,y_{11}\}$ with
\begin{eqnarray*}
y_1&=&2B_1-B_2+3B_6\\
y_4&=&-B_1+B_6\\
y_5&=&-2B_1+B_2-2B_6\\
y_7&=&-B_1-3B_2+B_5-12B_6\\
y_8&=&11B_1-B_3+B_5-11B_6\\
y_{11}&=&42B_1-21B_2-6B_3+B_4-6B_5+42B_6.
\end{eqnarray*}
For all integers $k$,
\begin{eqnarray*}
\psi^k(y_1)&=&ky_1-\tfrac12(k-k^4)y_4+\tfrac{11}{10}(k-k^5)y_5
+(\tfrac1{70}k-\tfrac{11}{120}k^5+\tfrac{13}{168}k^7)y_7\\
&&+(-\tfrac1{140}k+\tfrac1{480}k^4+\tfrac{11}{240}k^5-\tfrac{13}{336}
k^7-\tfrac1{480}k^8)y_8\\
&&+(\tfrac1{4620}k-\tfrac{11}{6720}k^5+\tfrac{13}
{13440}k^7+\tfrac{67}{147840}k^{11})y_{11}\\
\psi^k(y_4)&=&k^4y_4+\tfrac1{240}(k^4-k^8)y_8\\
\psi^k(y_5)&=&k^5y_5+\tfrac1{12}(k^5-k^7)y_7-\tfrac1{24}(k^5-k^7)y_8
+(\tfrac1{672}k^5-\tfrac1{960}k^7-\tfrac1{2240}k^{11})y_{11}\\
\psi^k(y_7)&=&k^7y_7-\tfrac12(k^7-k^8)y_8+\tfrac1{80}(k^7-k^{11})y_{11}\\
\psi^k(y_8)&=&k^8y_8\\
\psi^k(y_{11})&=&k^{11}y_{11}.
\end{eqnarray*}
\end{prop}

\begin{prop}\label{psiE7} A basis for $PK^1(E_7)$ is given by $\{y_1,y_5,
y_7,y_9,y_{11},y_{13},y_{17}\}$ with
\begin{eqnarray*} 
y_1&=&1873B_1-35B_2+15B_3-2B_5+287B_6-23056B_7\\
y_5&=&-113B_1-B_2+2B_3-B_5+29B_6-547B_7\\
y_7&=&6B_1-6B_2+7B_6-216B_7\\
y_9&=&2B_1-B_2+B_6-30B_7\\
y_{11}&=&-184B_1-B_2+3B_3-B_5+21B_6-292B_7\\
y_{13}&=&120B_1-5B_2-2B_3+B_5-22B_6+328B_7\\
y_{17}&=&1672B_1-252B_2-34B_3+B_4-12B_5+177B_6-1344B_7.
\end{eqnarray*}
For all integers $k$,
\begin{eqnarray*}\psi^k(y_1)&=&ky_1-\tfrac{22544}5(k-k^5)y_5+(\tfrac{29367}{70}k
-\tfrac{64814}{15}k^5+\tfrac{3659}{42}k^7)y_7+(\tfrac{35089}{21}k-\tfrac{57769}{45}k^5\\
&&-\tfrac{36590}{63}k^7+\tfrac{8714}{45}k^9)y_9+(-\tfrac{25409}{2310}k+\tfrac{1409}{210}k^5+\tfrac{3659}{1680}k^7+\tfrac{39031}
{18480}k^{11})y_{11}\\
&&+(-\tfrac{45223173}{10010}k+\tfrac{85358629}{18900}k^5
+\tfrac{464693}{181440}k^7+\tfrac{4357}{5400}k^9+\tfrac{507403}{221760}k^{11}
-\tfrac{24684827}{5896800}k^{13})y_{13}\\
&&+(\tfrac1{1021020}k+\tfrac{2154361}{154791000}k^5+\tfrac{3659}{684288}k^7
+\tfrac{4357}{1296000}k^9+\tfrac{39031}{11176704}k^{11}\\
&&-\tfrac{24684827}{707616000}k^{13}+\tfrac{8118218221}{926269344000}k^{17})y_{17}\\
\psi^k(y_5)&=&k^5y_5-\tfrac{23}{24}(k^5-k^7)y_7+(-\tfrac{41}{144}k^5
-\tfrac{115}{18}k^7+\tfrac{961}{144}k^9)y_9+(\tfrac1{672}k^5+\tfrac{23}{960}k^7\\
&&-\tfrac{57}{2240}k^{11})y_{11}+(\tfrac{60581}{60480}k^5+\tfrac{2921}{103680}k^7+\tfrac{961}{34560}k^9
-\tfrac{247}{8960}k^{11}-\tfrac{106799}{103680}k^{13})y_{13}\\
&&+(\tfrac{1529}{495331200}k^5+\tfrac{161}{2737152}k^7+\tfrac{961}{8294400}k^9
-\tfrac{19}{451584}k^{11}-\tfrac{106799}{12441600}k^{13}\\&&+\tfrac{1473021997}{174356582400}
k^{17})y_{17}\\
\psi^k(y_7)&=&k^7y_7-\tfrac{20}3(k^7-k^9)y_9+\tfrac1{40}(k^7-k^{11})y_{11}+(\tfrac{127}{4320}
k^7+\tfrac1{36}k^9-\tfrac{13}{480}k^{11}\\&&-\tfrac{13}{432}k^{13})y_{13}+(\tfrac7{114048}k^7+\tfrac1{8640}k^9-\tfrac1{24192}k^{11}-\tfrac{13}{51840}k^{13}
+\tfrac{17}{147840}k^{17})y_{17}\\
\psi^k(y_9)&=&k^9y_9+\tfrac1{240}(k^9-k^{13})y_{13}+(\tfrac1{57600}k^9-\tfrac1{28800}k^{13}
+\tfrac1{57600}k^{17})y_{17}\\
\psi^k(y_{11})&=&k^{11}y_{11}+\tfrac{13}{12}(k^{11}-k^{13})y_{13}+(\tfrac5{3024}k^{11}
-\tfrac{13}{1440}k^{13}+\tfrac{223}{30240}k^{17})y_{17}\\
\psi^k(y_{13})&=&k^{13}y_{13}+\tfrac1{120}(k^{13}-k^{17})y_{17}\\
\psi^k(y_{17})&=&k^{17}y_{17}
\end{eqnarray*}
\end{prop}

Before stating the final of these results about Adams operations in the
$K$-theory of exceptional Lie groups, the case in which many of the numbers become ridiculously large, we point out two features.
One is that the coefficient of each $y_j$ in $\psi^k(y_i)$ is actually an integer, yielding
integrality results. The other is that at least the second terms give information
about attaching maps via primes occurring in denominators. This follows
from \cite{JXIV}, where it was shown that in a sphere bundle over a sphere with attaching map $\a_t$ with $t\not\equiv0$ mod $p$ the Adams operations will be as described in our Proposition \ref{sphbdl}.
Thus, for example, the 5s in the denominators of the second terms in
the formulas for $\psi^k(y_1)$, $\psi^k(y_9)$, and $\psi^k(y_{13})$
in Proposition \ref{psiE7} are, at the very least, strongly suggestive
that there are $\a_1$ attaching maps from 1 to 5, 9 to 13, and 13 to 17 in $E_7$, and a similar deduction can be made at the prime 3.
The same conclusion can be made from somewhat simpler formulas of
$K(-)_{(p)}$, but here we are getting information about all primes at once. 

\begin{prop}\label{psiE8} A basis for $PK^1(E_8)$ is given by
$\{y_1,y_7,y_{11},y_{13},y_{17},y_{19},y_{23},y_{29}\}$ with
\begin{eqnarray*}
y_1&=&-784157B_1-30713B_2-218B_3+13B_5-919B_6+950380B_7-153687494B_8\\
y_7&=&224745B_1-2221B_2-101B_3+9B_5-950B_6+69688B_7-3701825B_8\\
y_{11}&=&100088B_1-982B_2-45B_3+4B_5-422B_6+30898B_7-1636120B_8\\
y_{13}&=&100091B_1-982B_2-45B_3+4B_5-422B_6+30897B_7-1635950B_8\\
y_{17}&=&71682B_1-788B_2-32B_3+3B_5-318B_6+23441B_7-1244530B_8\\
y_{19}&=&26482B_1-223B_2-12B_3+B_5-105B_6+7632B_7-403600B_8\\
y_{23}&=&28444B_1-195B_2-13B_3+B_5-104B_6+7462B_7-392340B_8\\
y_{29}&=&2691065B_1-38570B_2-1235B_3+B_4-58B_5+4409B_6\\
&&-174307B_7+6425670B_8.
\end{eqnarray*}
For all integers $k$,
\begin{eqnarray*}
\psi^k(y_1)&=&ky_1-\tfrac{1454165}{42}(k-k^7)y_7+(\tfrac{182723259}{154}k-\tfrac{650593421}{336}k^7+\tfrac{395881345}{528}k^{11})y_{11}\\
&&+(-\tfrac{185634674}{143}k+\tfrac{50198648299}{24192}k^7-
\tfrac{9105270935}{12672}k^{11}-\tfrac{9172441277}{157248}k^{13})y_{13}\\
&&+(\tfrac{1116883611}{6188}k-\tfrac{332043745267}{1596672}k^7
-\tfrac{4354694795}{145152}k^{11}+\tfrac{541174035343}{9434880}k^{13}\\
&&+\tfrac{14942613523}{135717120}k^{17})y_{17}+(\tfrac{10741728047}
{149226}k-\tfrac{113260168654603}{1634592960}k^7-\tfrac{326918814701}
{38320128}k^{11}\\
&&+\tfrac{1201589807287}{182891520}k^{13}-\tfrac{1240236922409}
{1628605440}k^{17}+\tfrac{2656402375049}{90348410880}k^{19})y_{19}\\
&&+(\tfrac{293377022893}{2028117}k-\tfrac{2311774903640460613}
{13338278553600}k^7-\tfrac{6031690056666301}{228248616960}k^{11}\\
&&+\tfrac{34318514761511237}{627683696640}k^{13}+\tfrac{354124997881577}
{820817141760}k^{17}-\tfrac{634880167636711}{21683618611200}k^{19}\\
&&+\tfrac{31696200305279479}{2454243253862400}k^{23})y_{23}+(\tfrac1
{38818159380}k+\tfrac{254146452881}{7362729761587200}k^7\\
&&+\tfrac{399916721902123}{72128388948295680}k^{11}-\tfrac{70600280509069}{10243797929164800}k^{13}+\tfrac{21980584492333}
{53779939128115200}k^{17}\\
&&+\tfrac{2656402375049}{5724475313356800}k^{19}+\tfrac{31696200305279479}{1236938599946649600}k^{23}-\tfrac{205974908943845743817}{8178638022847247155200}k^{27})y_{27}\\
\psi^k(y_7)&=&k^7y_7-\tfrac{2237}{40}(k^{7}-k^{11})y_{11}+(\tfrac{172603}{
2880}k^7-\tfrac{51451}{960}k^{11}-\tfrac{1825}{288}k^{13})y_{13}\\
&&+(-\tfrac{1141699}{190080}k^7-\tfrac{270677}{120960}k^{11}+\tfrac
{21535}{3456}k^{13}+\tfrac{535673}{266112}k^{17})y_{17}+(-\tfrac
{389433691}{194594400}k^7\\
&&-\tfrac{9236573}{14515200}k^{11}+\tfrac{621595}{870912}k^{13}
-\tfrac{44460859}
{3193344}k^{17}+\tfrac{8970785339}{566092800}k^{19})y_{19}\\
&&+(-\tfrac{7948805340661}{1587890304000}k^7-\tfrac{1874574277103}
{951035904000}k^{11}+\tfrac{1365640565}{229920768}k^{13}\\
&&+\tfrac{12694914427}{1609445376}k^{17}-\tfrac{2144017696021}
{135862272000}k^{19}+\tfrac{794109865426391}{88921857024000}k^{23})y_{23}\\
&&+(\tfrac{873857}{876515447808000}k^7+\tfrac{11299013179}{27321359450112000}k^{11}-\tfrac{561881}{750461386752}k^{13}\\
&&+\tfrac{787974983}{105450861035520}k^{17}+\tfrac{8970785339}
{35867639808000}k^{19}\\
&&+\tfrac{794109865426391}{44816615940096000}k^{23}
-\tfrac{168990789500875519}{9400733359594536960}k^{29})y_{29}\\
\psi^k(y_{11})&=&k^{11}y_{11}-\tfrac{23}{24}(k^{11}-k^{13})y_{13}
+(-\tfrac{121}{3024}k^{11}-\tfrac{1357}{1440}k^{13}+\tfrac{29707}
{30240}k^{17})y_{17}\\
&&+(-\tfrac{4129}{362880}k^{11}-\tfrac{39169}{362880}k^{13}
-\tfrac{2465681}{362880}k^{17}+\tfrac{2508979}{362880}k^{19})y_{19}\\
&&+(-\tfrac{837985819}{23775897600}k^{11}-\tfrac{86054063}{95800320}k
^{13}+\tfrac{704026193}{182891520}k^{17}-\tfrac{599645981}{87091200}
k^{19}\\
&&+\tfrac{129765607619}{32691859200}k^{23})y_{23}+(\tfrac{5050967}
{683033986252800}k^{11}+\tfrac{177031}{1563461222400}k^{13}\\
&&+\tfrac{43698997}{11983052390400}k^{17}+\tfrac{228089}{2090188800}k^{1
9}+\tfrac{129765607619}{16476697036800}k^{23}-\tfrac{371040821209573}
{46446311065190400}k^{29})y_{29}\\
\psi^k(y_{13})&=&k^{13}y_{13}-\tfrac{59}{60}(k^{13}-k^{17})y_{17}
+(-\tfrac{1703}{15120}k^{13}-\tfrac{4897}{720}k^{17}+\tfrac{5227}{756}
k^{19})y_{19}\\
&&+(-\tfrac{3741481}{3991680}k^{13}+\tfrac{1398241}{362880}k^{17}
-\tfrac{1249253}{181440}k^{19}+\tfrac{3961099}{997920}k^{23})y_{23}\\
&&+(\tfrac{7697}{65144217600}k^{13}+\tfrac{86789}{23775897600}k^{17}
+\tfrac{5227}{47900160}k^{19}\\
&&+\tfrac{3961099}{502951680}k^{23}-\tfrac{142072160653}{17784371404800}k^{29})y_{29}\\
\psi^k(y_{17})&=&k^{17}y_{17}-\tfrac{83}{12}(k^{17}-k^{19})y_{19}
+(\tfrac{23699}{6048}k^{17}-\tfrac{19837}{2880}k^{19}+\tfrac{179587}
{60480}k^{23})y_{23}\\
&&+(\tfrac{1471}{396264960}k^{17}+\tfrac{83}{760320}k^{19}
+\tfrac{179587}{30481920}k^{23}-\tfrac{26172961}{4358914560}k^{29})y_{29}\\
\psi^k(y_{19})&=&k^{19}y_{19}-\tfrac{239}{240}(k^{19}-k^{23})
+(\tfrac1{63360}k^{19}+\tfrac{239}{120960}k^{23}-\tfrac{265}{133056}
k^{29})y_{29}\\
\psi^k(y_{23})&=&k^{23}y_{23}+\tfrac1{504}(k^{23}-k^{29})y_{29}\\
\psi^k(y_{29})&=&k^{29}y_{29}
\end{eqnarray*}
\end{prop}

\section{Determination of $\vp_{2m}(E_8;5)$}\label{even}
Bousfield proved the following result in \cite{Bous}. This is the
result that has been called \lq\lq Bousfield's theorem'' throughout this paper.
\begin{thm}\label{Bousthm} Let $X$ be a 1-connected finite $H$-space with
$H_*(X;\bq)$ associative, $p$ an odd prime, and $r$ a generator of
$(\bz/p^2 )^\times$. Then
$$\vp_{2m}(X;p)=(\coker(\phi_m))^\#\text{ and }\vp_{2m-
1}(X;p)=(\ker(\phi_m))^\#,$$
where $(-)^\#$ denotes the Pontryagin dual, and
$$\phi_m=\psi^r-r^m:PK^1(X;\zphat)/\im(\psi^p)\to PK^1(X;\zphat)/\im(\psi^p).$$
\end{thm}
For the finite abelian groups with which we deal, $PK^1(X;\zphat)/\im(\psi^p)$
is isomorphic to $PK^1(X;\bz_{(p)})/\im(\psi^p)$, and the only effect of the
Pontryagin dual
is to reverse the direction of arrows.

The following result is immediate from Proposition \ref{niceops} and Theorem
\ref{Bousthm}.
\begin{prop} $\vp_{2m}(E_8;5)\approx A_m\oplus B_m$, where
$$A_m=\bz_{(5)}(x_1,x_{13},x_{17},x_{29})/(r_1,\ldots,r_8)$$
and $$B_m=\bz_{(5)}(x_7,x_{11},x_{19},x_{23})/(s_1,\ldots,s_8),$$
where $\bz_{(5)}(-)$ denotes the free $\bz_{(5)}$-module on the indicated
generators, and the relations are given as follows:
\begin{eqnarray*}
r_1&:& 5^{29}x_{29};\\
r_2&:& 5^{17}x_{17}-5^{16}(5^{12}-1)x_{29};\\
r_3&:& 5^{13}x_{13}+4\cdot 5^{12}(5^4-1)x_{17}+2\cdot 5^{11}(1+2\cdot5^4-3\cdot
5^{16})x_{29};\\
r_4&:& 5x_1+10\cdot5(5^{12}-1)x_{13}-8\cdot5(5^{12}-1)x_{17}+4(5^{12}-
1)x_{29};\\
r_5&:& (2^{29}-2^m)x_{29};\\
r_6&:& (2^{17}-2^m)x_{17}-\tfrac15\cdot2^{17}(2^{12}-1)x_{29};\\
r_7&:& (2^{13}-2^m)x_{13}+12\cdot 2^{13}x_{17}+\tfrac3{25}\cdot2^{14}(11-
2^{16})x_{29};\\
r_8&:& (2-2^m)x_1+20(2^{12}-1)x_{13}-16(2^{12}-1)x_{17}+\tfrac85(2^{12}-
1)x_{29};\\
s_1&:& 5^{23}x_{23};\\
s_2&:& 5^{19}x_{19}+5^{18}(1-5^4)x_{23};\\
s_3&:& 5^{11}x_{11}+5^{10}(1-5^8)x_{19}+5^9(9\cdot 5^{12}-5^8-8)x_{23};\\
s_4&:& 5^7x_7+2\cdot5^6(1-5^4)x_{11}+5^5(3-2\cdot5^4-5^{12})x_{19}\\
&&+5^4(32+16\cdot 5^4-5^{12}-47\cdot5^{16})x_{23};\\
s_5&:&(2^{23}-2^m)x_{23};\\
s_6&:&(2^{19}-2^m)x_{19}-3\cdot2^{19}x_{23};\\
s_7&:&(2^{11}-2^m)x_{11}-51\cdot2^{11}x_{19}+183\cdot2^{14}x_{23};\\
s_8&:&(2^7-2^m)x_7-3\cdot2^8x_{11}-165\cdot2^7x_{19}-771\cdot2^{12}x_{23}.
\end{eqnarray*}
\label{16r}\end{prop}

We shall analyze $A_m$ first. We will make frequent use of the following
well-known fact proved in \cite[2.12]{JX}.
\begin{prop}\label{nu}
If $r$ generates the group of units $(\bz/p^2)^\times$, then
$$\nu_p(r^m-1)=\begin{cases}0&\text{if $m\not\equiv 0$ mod $p-1$}\\
1+\nu_p(m)&\text{if $m\equiv0$ mod $p-1$}.\end{cases}$$
\end{prop}

Relations $r_5$, $r_6$, $r_7$, and $r_8$ imply immediately that $A_m$ is 0 if
$m\not\equiv1$ mod 4. We will write $m=4k+1$ and divide $r_5$, $r_6$, $r_7$,
and $r_8$ by the units $2^{29}$, $2^{17}$, $2^{13}$, and $2$, respectively.
We use $r_6$ to eliminate $x_{29}$, and replace it by $((1-2^{4k-
16})/819)x_{17}$ in the other relations, and then use the modified $r_7$ to
eliminate $x_{17}$, replacing it by a specific multiple of $x_{13}$ in all
relations. Next we add $\frac15(2^{4k}-1)r_4$ to $r_8$ to eliminate $x_1$ from
$r_8$. Now the modified $r_4$ is the only relation involving $x_1$, and it is
of the form $5x_1+\a5^2x_{13}$. This implies that $x_1+5\a x_{13}$ generates
a $\bz/5$ direct summand. The remaining relations on $x_{13}$ are as follows,
obtained respectively from $r_1$ and $r_5$, $r_2$, $r_3$, and $r_8$.
\begin{eqnarray*}
t_1&:&5^{\min(28,\nu(k-7))+\nu(k-4)+\nu(k-3)+3}\\
t_2&:&5^{16}(4095+(5^{12}-1)(2^{4(k-4)}-1))\\
t_3&:&5^{13}\biggl(6\cdot819+3\cdot2621(2^{4k-16}-1)\\
&&\quad +\tfrac15(2^{4k-12}-1)
\bigl(2(5^4-1)819-\tfrac15(2^{4k-16}-1)(1+2\cdot5^4-
3\cdot5^{16})\bigr)\biggr)\\
t_4&:&5^2
\biggl(12\cdot819+6\cdot2621(2^{4k-16}-1)-\tfrac15(2^{4k-12}-1)(4\cdot819+
\tfrac25(2^{4k-16}-1))\\
&&\quad +\tfrac15(2^{4k}-1)(5^{12}-1)
\biggl(12+\tfrac{6\cdot2621}{819}(2^{4k-16}-1)\\
&&\quad-\tfrac25(2^{4k-12}-1)(2+\tfrac15(2^{4k-16}-1)/819)\biggr)\biggr)
\end{eqnarray*}

The relation $t_4$ can be manipulated to a unit times
\begin{equation}\label{t4n}\tfrac15(2^{4k-16}-1)(2^{4k-28}-1)(2^{12}
(2^{4k-12}-1)-5^{12}(2^{4k}-1)),\end{equation}
while $t_3$ can be rewritten as
\begin{equation}\label{t3n}-2^{16}5^{11}(2^{4k-16}-1)(2^{4k-28}-1)
-2^{13}5^{15}(2^{4k-12}-1)(2^{4k-28}-1)+3\cdot5^{28}(2^{4k-16}-
1),\end{equation}
and $t_2$ as
\begin{equation}\label{t2}5^{28}(2^{4k-16}-1)-5^{16}2^{12}(2^{4k-28}-1).
\end{equation}
These four relations on $x_{13}$ result in a summand of order $5^e$, where
\begin{itemize}
\item if $k\equiv0,1$ mod 5, then $e=2$ from (\ref{t4n});
\item if $k\equiv3$ mod 5, then $e$ is the minimum of
 $2+\nu(k-3)$ from (\ref{t4n}) or $13$ from
(\ref{t3n});
\item if $k\equiv4$ mod 5, then $e$ is the minimum of $2+\nu(k-4)$
from (\ref{t4n}) or $17$ from (\ref{t2});
\item if $k\equiv2$ mod 5, then $e$ is the minimum of $2+\nu(k-7)$
from (\ref{t4n}) or $29$ from (\ref{t3n}).
\end{itemize}
Thus $A_m$ gives the groups claimed in Theorem \ref{main} in the case $m\equiv
1$ mod 4.

Now we use the relations $s_1,\ldots,s_8$ to determine $B_m$. The analysis is
similar to that just employed
for $A_m$, but we use a different technique near the
end. First we use $s_5,\ldots,s_8$ to see that $B_m=0$ if $m\not\equiv3$ mod 4.
Now we write $m=4k+3$, and divide the last four relations by their initial
2-power. We use $s_6$ to eliminate $x_{23}$, and then $s_7$ to eliminate
$x_{19}$. When these expressions are substituted into $s_8$, it becomes
$$(1-2^{4k-4})x_7-ux_{11},$$
where
\begin{equation}\label{u}
u=6+\bigl(165+257\cdot2^5(1-2^{4k-16})\bigr)(1-2^{4k-8})/
(488\cdot 2^{4k-16}-437)\end{equation}
is a unit.
Thus $x_{11}$ can be eliminated, and so $B_m$ is cyclic with generator
$x_7$ and relations
\begin{eqnarray*}
w_1&:&5^{\min(26,4+\nu(k-5)+\nu(k-4)+\nu(k-2)+\nu(k-1))};\\
w_2&:&5^{20+\nu(k-2)+\nu(k-1)}(15+(1-5^4)(1-2^{4k-16}));\\
w_3&:&5^{10+\nu(k-1)}\biggl(5^2(488\cdot2^{4k-16}-437)\\
&&\quad +(1-2^{4k-8})\bigl(5-5^9+\tfrac13(1-2^{4k-16})(9\cdot5^{12}-5^8-8)
\bigr)\biggr);\\
w_4&:&5^4\biggl
(5^3u+(1-2^{4k-4})\biggl(2(5^2-5^6)+\frac{(1-2^{4k-8})}{488\cdot2^{4k-16}-
437}\\
&&\quad
\cdot\bigl(15-2\cdot5^5-5^{13}+\tfrac13(1-2^{4k-16})(32+16\cdot5^4-5^{12}-
47\cdot5^{16})\bigr)\biggr)\biggr),
\end{eqnarray*}
where $u$ is as in (\ref{u}).
We easily read off from these that the order of $B_m$ is $5^{\min(7,4+\nu(k-
1))}$ if $k\equiv 1$ or 3 mod 5.

Now let $P=2^{4k}$, and write $w_4$ as a unit times $5^4$ times the following
expression.
\begin{eqnarray*}
&&5^3\biggl(6\cdot2^{12}(488P-437\cdot2^{16})+2^4\bigl
(165\cdot2^{16}+257\cdot32(2^{16}-
P)\bigr)(2^8-P)\biggr)\\
&&\quad +(16-P)
\biggl(50\cdot2^8(1-5^4)(488P-437\cdot2^{16})+(2^8-P)\bigl(2^{16}(15-
2\cdot5^5-5^{13})\\
&&\quad+\tfrac13(2^{16}-P)(32+16\cdot5^4-5^{12}-47\cdot5^{16})\bigr)\biggr).
\end{eqnarray*}
We use {\tt Maple}
to write this expression as $A+B\cdot P+C\cdot P^2+D\cdot P^3$,
where $A$, $B$, $C$, and $D$ are certain explicit large integers. We note that
for any positive integer $e$, this cubic expression can be rewritten as
$$a_0+a_1(P-2^e)+a_2(P-2^e)^2+a_3(P-2^e)^3,$$
where
\begin{eqnarray*}
a_3&=&D\\
a_2&=&C+3\cdot2^eD\\
a_1&=&B+2^{e+1}C+3\cdot2^{2e}D\\
a_0&=&A+2^eB+2^{2e}C+2^{3e}D.
\end{eqnarray*}
With $e=8$, {\tt Maple}
computes these $a_i$, from which we deduce that $w_4$ can be
written as
$$5^4(u_15^7+u_25^3(2^{4k}-2^8)+u_35(2^{4k}-2^8)^2+u_4(2^{4k}-2^8)^3),$$
where $u_i$ are units in $\bz_{(5)}$.
The relation given by this and $w_1$ when $k\equiv2$ mod 5 is
the desired $5^{\min(11,4+\nu(k-2))}$, and $w_2$ and $w_3$ are easily seen
to imply no additional restrictions.

If $k\equiv4$ mod 5, we use $e=16$ in the above method (with the same values
of $A$, $B$, $C$, and $D$) and obtain the following form of $w_4$, where again
$u_i$ are units in $\bz_{(5)}$.
$$5^4\bigl(u_15^{15}+u_25^2(2^{4k}-2^{16})+u_35^2(2^{4k}-2^{16})^2+u_4(2^{4k}-
2^{16})^3\bigr)$$
This and $w_1$ give $5^{\min(19,4+\nu(k-4))}$ as the relation. To see that
$w_3$ can give no smaller relation, we can write it as
$$5^{10}((2^{4k-16}-1)(u_5(2^{4k-16}-1)+5u_6)+5^{10}u_7),$$
with $u_i$ units.

If $k\equiv5$ mod 5, we use $e=20$ in the above method and find that $w_4$ can
be written as
$$5^4(u_15^{19}+u_25^2(2^{4k}-2^{20})+u_35(2^{2k}-2^{20})^2+u_4(2^{2k}-
2^{20})^3).$$
This and $w_1$ give the desired $5^{\min(24,4+\nu(k-5))}$. To see that $w_2$
cannot give a smaller relation in this case, we write it as a unit times
$$5^{20}(2^4(2^{4k-20}-1)+5^4(2^{4k-16}-1)).$$
Finally, $w_3$ may be rewritten, using the method employed above for $w_4$, as
$$5^{24}u_0+5^{11}u_1(2^{2k}-2^{20})+5^{10}u_2(2^{4k}-2^{20})^2,$$
which implies that it cannot give a smaller relation here.

This completes the calculation of $B_m$, showing that it yields exactly the
values of $\vp_{2m}(E_8;5)$ claimed in Theorem \ref{main} when $m\equiv 3$ mod
4. This completes the determination of $\vp_*(E_8;5)$ when $*$ is even.

\section{Determination of $\vp_{2m-1}(E_8;5)$}\label{odd}
We could compute
$\vp_{2m-1}(E_8;5)$ using Theorem \ref{Bousthm} by computing the
relevant kernel. These computations are probably even more difficult than those
of Section \ref{even}, and so, uncharacteristically, we follow a program of
avoiding calculations as much as possible. A sequence of five propositions will
yield the result for $\vp_{2m-1}(E_8;5)$ claimed in Theorem \ref{main}.

The first proposition is elementary and well-known.
\begin{prop}\label{order} Let $p$ be odd. If $X$ is a simply-connected
finite $H$-space with $H_*(X;\bq)$ associative or $X$ is a space built from a
finite number of odd-dimensional spheres by fibrations, then
$\vp_{2m-1}(X;p)$ and $\vp_{2m}(X;p)$ have the same order.\end{prop}
\begin{pf} For the first type of $X$, we use Theorem \ref{Bousthm} and the
observation that the kernel and cokernel of an endomorphism of a finite
abelian group have the same order. For the second type of $X$, we use the
argument employed for $SU(n)$ in \cite[p.530]{D}.\end{pf}

Next we have the following recent observation of Bendersky.
\begin{prop}\label{number} If $X$ is a space with an $H$-space exponent at $p$,
the number of direct summands in $\vp_i(X;p)$ depends only on the residue mod
$q$ of $i$.\end{prop}
Here we begin
using $q=2(p-1)$, and say that $X$ has an $H$-space exponent at $p$ if
some $p^e$-power
map on some iterated loop space of $X$ is null-homotopic. Compact
Lie groups and spheres have $H$-space exponents at all primes.
\begin{pf*}{Proof of Proposition \ref{number}.}
There is a commutative diagram of isomorphisms
$$\begin{matrix} \pi_i(X;\bz/p)&\approx&\pi_i(X)\ot\bz/p\oplus\tor(\pi_{i-
1}(X),\bz/p)\\
\da&&\da\\
\pi_{i+kqp^e}(X;\bz/p)&\approx&\pi_{i+kqp^e}(X)\ot\bz/p\oplus\tor(\pi_{i-
1+kqp^e}(X),\bz/p)
\end{matrix}$$
where the vertical arrows are as described in \cite[\S2]{D}. These arrows are
defined using the null-homotopy of the $p^e$-power map.
The one on the left can equivalently be defined using Adams maps of the mod $p$
Moore space. The direct limit of these vertical morphisms
defines the $v_1$-periodic homotopy groups. The diagram commutes because the
morphisms $\pi_*X\to\pi_{*+kqp^e}X$ commute with multiplication by $p$.

Passing to the direct limit, we obtain an isomorphism
$$\vp_i(X;\bz/p)\approx \vp_i(X)\ot\bz/p\oplus\tor(\vp_{i-1}(X),\bz/p).$$
These will be finite $\bz_{(p)}$-modules, and if we denote by $\#(G)$
the number of direct summands in such a module, we obtain
$$\#(\vp_i(X;\bz/p))=\#(\vp_i(X))+\#(\vp_{i-1}(X)).$$

On the other  hand, $\vp_i(X;\bz/p)\approx\vp_{i+q}(X;\bz/p)$ because the
$v_1^{p^e}$-map which defines the direct limit can be chosen to be the
composite of maps $\Sigma^qM(p)\to M(p)$. Thus $\#(\vp_i(X))+\#(\vp_{i-1}(X))$
has period $q$. This implies that $\#(\vp_i(X))$ has period $q$. To see this, we note that if
$\Delta=\#(\vp_{q}X)-\#(\vp_0X)$, it follows that
$$\#(\vp_{j+q}X)-\#(\vp_jX)=(-1)^{j}\Delta$$
for all $j$, which implies that
$$\#(\vp_{j+Lq}(X))-\#(\vp_{j}(X))=(-1)^jL\Delta.$$
If $\Delta\ne0$, let $j=1$ if $\Delta>0$ and $j=0$ if $\Delta<0$. This implies the absurd statement
that $\#(\vp_{j+Lq}(X))<0$ for $L$ sufficiently large, and so we conclude $\Delta=0$.
\end{pf*}

By these two propositions and the computation of Section \ref{even},
Theorem \ref{main} will be proven once we show $\vp_{2m-1}(E_8;5)$ is cyclic
for some $m\equiv1$ mod 4 and for some $m\equiv3$ mod 4. This could be done by
an explicit calculation, but we prefer to approach it by exact sequences.
The following result is immediate from the Snake Lemma.
\begin{prop} Suppose $0\to M_1\to M_2\to M_3\to 0$ is a short exact sequence of
$p$-adic Adams modules with $\psi^p$ acting injectively on each. Let
$$\phi_i=\psi^r-r^m:M_i/\im(\psi^p)\to M_i/\im(\psi^p).$$
Then there is an exact sequence
$$0\to\ker\phi_1\to\ker\phi_2\to\ker\phi_3\to\coker\phi_1\to\coker\phi_2\to
\coker\phi_3\to0.$$
\label{psiles}\end{prop}
Theorem \ref{Bousthm} and exactness of Pontryagin duality yield the following
corollary.
\begin{cor} Suppose
maps $X\to Y\to Z$ induce a short exact sequence of Adams
modules
$$0\to PK^1(Z;\zphat)\to PK^1(Y;\zphat)\to PK^1(X;\zphat)\to0,$$
with $\psi^p$ acting injectively on each. Then there is a long exact sequence
$$\to\vp_{i+1}Z\to \vp_iX\to \vp_iY\to \vp_iZ\to \vp_{i-1}X\to.$$
\label{les}\end{cor}

Next we have the following basic calculation, which is the Bousfield approach
to a sphere bundle over sphere with attaching map $\a_t$. This result is
analogous to \cite[1.3]{BDMi}, which was obtained using the UNSS.
\begin{prop}\label{sphbdl} If $PK^1(X;\zphat)$ has generators $x$ and $y$ with
$\psi^ky=k^{n+t(p-1)}y$ and $\psi^kx=k^nx+\frac upk^n(k^{t(p-1)}-1)y$, where
$u$ is a unit in $\bz_{(p)}$, and $t\not\equiv0$ mod $p$,
then the only nonzero groups $\vp_*(X;p)$ are
$$\vp_{2n+qi}(X;p)\approx \vp_{2n+qi-1}(X;p)\approx\bz/p^e,$$
where
$$e=\begin{cases}\min(\nu_p(i)+2,n)&\text{if $i\not\equiv t$ mod $p$}\\
\min(\nu_p(i-t)+2,t(p-1)+n)&\text{if $i\equiv t$ mod $p$, and $n>1$}\\
\min(\nu_p(i-t-tp^{t(p-1)})+1,t(p-1)+2)&\text{if $i\equiv t$ mod $p$, and $n=1$}
\end{cases}$$
\end{prop}
\begin{pf} We again make frequent use of Proposition \ref{nu}.
By Theorem \ref{Bousthm}, $\vp_{2m}(X;p)$ has generators $x$ and $y$ subject to
relations
\begin{eqnarray*}
&&p^{n+t(p-1)}y\\
&&p^nx+up^{n-1}(p^{t(p-1)}-1)y\\
&&(r^{n+t(p-1)}-r^m)y\\
&&(r^n-r^m)x+\tfrac upr^n(r^{t(p-1)}-1)y
\end{eqnarray*}
The last two relations imply that the group is 0 unless $m\equiv n$ mod $p-1$.
We let $m=n+(p-1)i$. Since $t\not\equiv0$ mod $p$,
the fourth relation can be used to eliminate $y$. We
obtain that the group is cyclic with generator $x$, and relations on $x$
\begin{eqnarray*}&&p^{\nu(i)+1+\min(\nu(i-t)+1,n+t(p-1))}\\
&&(r^{t(p-1)}-1)p^{n-1}-p^{n-1}(p^{t(p-1)}-1)(1-r^{(p-1)i}).
\end{eqnarray*}
The second of these relations can be rewritten as
$$p^{n-1}\bigl(r^{t(p-1)}(1-r^{(p-1)(i-t)})+p^{t(p-1)}(r^{(p-1)t}-1)\bigr).$$
Now inspection of the relations shows that if $i-t\not\equiv0$ mod $p$, the
first relation becomes $p^{\nu(i)+2}$, and the second relation becomes
$p^n$. If $i-t\equiv0$ mod $p$, then the first relation becomes
$p^{\min(\nu(i-t)+2,n+t(p-1)+1)}$, while the second becomes
$p^{n+\min(\nu(i-t),t(p-1))}$ unless $\nu(i-t)=t(p-1)$.
If $i-t=\a p^{t(p-1)}$ and we write $r^{p-1}=1+\b p$ with $\b\not\equiv0$
mod $p$, then the second relation becomes $(i-\a)\b p^{t(p-1)+1}$ mod
$p^{t(p-1)+2}$. It is now routine to translate these relations when
$i-t\equiv0$ mod $p$ into those claimed in the proposition. This establishes
the results for the even-dimensional groups.

The odd groups could be computed directly; however, we can avoid computation as
follows. By Proposition \ref{order}, the odd groups have the asserted order,
and by Proposition \ref{number} it suffices to show that any one of them is
cyclic. By Proposition \ref{psiles}, there is an exact sequence
\begin{eqnarray}
&&0\to\vp_{2m}S^{2n+1}\to\vp_{2m}X\to\vp_{2m}S^{2n+1+tq}\mapright{\delta}
 \vp_{2m-1}S^{2n+1}\nonumber\\
&&\qq\to \vp_{2m-1}X\to\vp_{2m-1}S^{2n+1+tq}\to0,\label{Xseq}
\end{eqnarray}
where $\vp_*S^{2m+1}$ are the groups associated to the Adams module
$PK^1(S^{2m+1};\zphat)$, which is cyclic with $\psi^k=\cdot k^m$. It is easy to
use Theorem \ref{Bousthm} to compute that the only nonzero groups
$\vp_*S^{2m+1}$ are
\begin{equation}\vp_{2m+iq}S^{2m+1}\approx\vp_{2m+iq-
1}S^{2m+1}\approx\bz/p^{\min(m,\nu(i)+1)},\label{sph}\end{equation}
in agreement with well-known results. Choose $m=n+(p-1)(t+p^{n+t(p-1)})$.
By (\ref{sph}) and our computation of $\vp_{\text{ev}}X$, the exact sequence
(\ref{Xseq}) becomes
\begin{eqnarray*}
&&0\to\bz/p\to\bz/p^{n+(p-1)t}\to\bz/p^{n+(p-1)t}\mapright{\delta}
\bz/p\\
&&\to\vp_{2n+qt+qp^{n+t(p-1)}-1}X\to\bz/p^{n+(p-1)t}\to0,\end{eqnarray*}
which implies that $\delta$ is surjective and $\vp_{2n+qt+qp^{n+t(p-1)}-1}X
\approx\bz/p^{n+(p-1)t}$.
\end{pf}

Now, the determination of $\vp_{\text{odd}}(E_8;5)$ is completed by
using the observation which precedes Proposition \ref{psiles} and the following
result.
\begin{prop} If $m=4k+1$ and $k\equiv0$ or $1$ mod $5$, then $\vp_{2m-
1}(E_8;5)\approx\bz/5^3$. If $m=4k+3$ and $k\equiv3$ mod $5$, then
$\vp_{2m-1}(E_8;5)\approx\bz/5^4$.\label{twoc}\end{prop}
\begin{pf} Let $M_1^{29}$ (resp. $M_7^{23}$) denote the sub Adams module of
$PK^1(E_8)_{(5)}$ generated by $x_1$, $x_{13}$, $x_{17}$, and $x_{29}$
(resp. $x_7$, $x_{11}$, $x_{19}$, and $x_{23}$). (See Proposition
\ref{niceops}.) If $i,j\in\{1,13,17,29\}$ with $i\le j$, let $M_i^j$ denote the
subquotient of $M_1^{29}$ generated by those $x_k$ in $M_1^{29}$ with
$i\le k\le j$, and similarly for subquotients of $M_7^{23}$. Let
$$\phi_i^j=\psi^2-2^m:M_i^j/\im(\psi^5)\to M_i^j/\im(\psi^5)$$
for a fixed value of $m$ of the type specified in the proposition.
Let $K_i^j=\ker(\phi_i^j)$ and $C_i^j=\coker(\phi_i^j)$.

We consider first the case $m=4k+1$ with $k\equiv 0$ or 1 mod $5$.
By (\ref{sph}) we have $K_i^i\approx C_i^i\approx\bz/5$ for
$i\in\{1,13,17,29\}$, and by Propositions \ref{sphbdl} and \ref{niceops},
$$K_{13}^{17}\approx C_{13}^{17}\approx K_{17}^{29}\approx C_{17}^{29}\approx
\bz/5^2.$$
By Theorem \ref{main}, $C_1^{29}\approx \bz/5\oplus\bz/5^2$, and, by a
calculation similar to that of $A_m$ in Section \ref{even} but much easier,
we obtain
\begin{equation}\label{125} C_{13}^{29}\approx \bz/5^3.\end{equation}
We will sketch this calculation at the end of this section. The results about
$C$- and $K$-groups just listed imply, by just diagram chasing, that
\begin{equation}\label{129}K_1^{29}\approx \bz/5^3,
\end{equation} which implies the first part of the proposition
by Theorem \ref{Bousthm}.

Before we present the simple proof of (\ref{129}), we wish to present the
motivation or underlying rationale for this argument and these results.
It involves homotopy charts of the type used extensively in \cite{E7}.

\begin{diag}\label{firstdiag}
\begin{center}
\begin{picture}(365,65)
\def\mp{\multiput}
\def\elt{\circle*{3}}
\put(150,10){$59$}
\put(150,25){$35$}
\put(150,40){$27$}
\put(155,55){$3$}
\put(163,0){$\begin{scriptsize}2m-1\end{scriptsize}$}
\put(205,0){$\begin{scriptsize}2m\end{scriptsize}$}
\mp(180,15)(0,15){4}{\elt}
\mp(210,15)(0,15){4}{\elt}
\mp(180,15)(30,0){2}{\line(0,1){30}}
\put(210,15){\vector(-2,3){30}}
\end{picture}
\end{center}
\end{diag}

For the spheres $S$ of dimensions 3, 27, 35, and 59 whose $PK^1$-Adams modules
build this portion of $PK^1(E_8)_{(5)}$, and for the value of $m$ being
considered here,
$$\vp_{2m}(S)\approx\vp_{2m-1}(S)\approx\bz/5.$$
These are indicated by the dots. The vertical lines from 59 to 35, and from 35
to 27, are nontrivial multiplication by 5 obtained from Proposition
\ref{sphbdl}, which applies because of the terms $-\frac15k^{17}(k^{12}-
1)x_{29}$ and
$\frac45k^{13}(k^4-1)x_{17}$ in Proposition \ref{niceops}. (The first
of these is the all-important term which implies that Harper's product
decomposition was false.) The computation in Section \ref{even}
of $A_m\approx \bz/5\oplus\bz/5^2$ (for appropriate $m$) implies that there
must be a differential from $\vp_{2m}(S^{59})$, but it wasn't clear whether it
hit $\vp_{2m-1}(S^{27})$ or $\vp_{2m-1}(S^3)$. (Actually, the UNSS or the
$\frac45k(k^{12}-1)$-term in \ref{niceops} make it pretty clear that it hits
$\vp_{2m-1}(S^3)$.) Our computation that $C_{13}^{29}\approx\bz/5^3$ implies
that there was no differential from 59 to 27, and so it must go from 59 to 3,
as indicated in the diagram by the diagonal line.

But this is just one way of thinking. The result (\ref{129}) can be obtained
by diagram-chasing as follows. The claims made above about certain
$K$-groups and $C$-groups imply that the relevant exact sequences of
\ref{psiles} yield short exact sequences in the following commutative diagram,
with all groups except the middle one given explicitly in the second diagram.

$$\begin{CD} @.0@.0@.0@.\\
@.@VVV@VVV@VVV@.\\
0@>>> K_{29}^{29}@>>> K_{17}^{29}@>>> K_{17}^{17}@>>>0\\
@.@V\approx VV@VVV@VVV@.\\
0@>>>K_{29}^{29}@>>> K_{13}^{29}@>>>K_{13}^{17}@>>>0\\
@.@.@VVV@VVV@.\\
@.@.K_{13}^{13}@>\approx >>K_{13}^{13}\\
@.@.@VVV@VVV@.\\
@.@.0@.0@.
\end{CD}$$
\bigskip
\bigskip
$$\begin{CD}\bz/5@>>> \bz/5^2@>>> \bz/5\\
@V\approx VV@VVV@VVV\\
\bz/5@>>>K_{13}^{29}@>>>\bz/5^2\\
@.@VVV@VVV\\
@.\bz/5@>\approx >>\bz/5
\end{CD}$$

It is easy to verify that in such a diagram, we must have $K_{13}^{29}\approx
\bz/5^3$. Another exact sequence from \ref{psiles} begins
$$0\to K_{13}^{29}\to K_1^{29}\to.$$
Since $|K_1^{29}|=|C_1^{29}|=5^3$, and $K_1^{29}$ contains a cyclic subgroup,
$K_{13}^{29}$, of
order $5^3$, we must have $K_1^{29}\approx\bz/5^3$, as claimed in (\ref{129}).

The computation when $m=4k+3$ and $k\equiv3$ mod 5 is easier.  In this case,
each of the four spheres yields a $\bz/5$, and each pair of consecutive spheres
has a nontrivial $\cdot 5$ extension because of \ref{sphbdl} and the terms
$\frac25k^7(1-k^4)x_{11}$, $\frac15k^{11}(1-k^8)x_{19}$, and $\frac15k^{19}(1-
k^4)x_{23}$ in \ref{niceops}. Since $|K_7^{23}|=|C_7^{23}|=5^4$, we can deduce
cyclicity of $K_7^{23}$ by chasing diagrams such as those used just above.
The relevant homotopy chart is as below, with no differential.

\begin{center}
\begin{picture}(65,65)
\def\mp{\multiput}
\def\elt{\circle*{3}}
\put(0,10){$47$}
\put(0,25){$39$}
\put(0,40){$23$}
\put(0,55){$15$}
\put(13,0){$\begin{scriptsize}2m-1\end{scriptsize}$}
\put(55,0){$\begin{scriptsize}2m\end{scriptsize}$}
\mp(30,15)(0,15){4}{\elt}
\mp(60,15)(0,15){4}{\elt}
\mp(30,15)(30,0){2}{\line(0,1){45}}
\end{picture}
\end{center}

We close this section by sketching the proof of (\ref{125}). We are computing
$A_m$ in Proposition \ref{16r}, with $x_1$, $r_4$, and $r_8$ removed, and
$m=4k+1$ with $k\equiv0$ or 1 mod 5. In the analysis in Section \ref{even},
we will have generator $x_{13}$ with relations $t_1$, $t_2$, and $t_3$.
The relation $t_1$ says that $5^3x_{13}=0$, and the other relations involve
much larger powers of 5. Thus the group is $\bz/5^3$, as claimed.
\end{pf}

\section{Calculation of $\vp_*(E_8;3)$}\label{3calc}
The computation of $\vp_*(E_8;3)$ is performed similarly to that of
$\vp_*(E_8;5)$. One thing that makes the analysis more complicated is that all eight generators are related to one another by Adams operations, rather than being divided into two groups of four, as was the case for $(E_8,5)$. It is more difficult to analyze an abelian group with 8 generators and 16 relations than to do it for two groups, each with 4 generators and 8 relations. We rely on {\tt Maple} for every step of this computation.

For $i\in\{1,7,11,13,17,19,23,29\}$, let $v_i$ be the columns of (\ref{evecs}), satisfying $\psi^k(v_i)=k^iv_i$. Then 32 times we replace vectors $v$ by $v':=(v-w)/3$, with $w$ a linear combination of vectors in the set which follow $v$, and $v'$ still integral, and obtain finally a basis $\{w_1,w_7,w_{11},w_{13},w_{17},w_{19},w_{23},w_{29}\}$ for
$PK^1(E_8;\bz_{(3)})$ defined by
\begin{eqnarray*}
w_1&=&v_1/3^2+730v_7/3^8+10v_{11}/3^{10}+4v_{13}/3^8+28v_{17}/3^9-109v_{19}/3^8\\
&&+1364v_{23}/3^9-44078v_{29}/3^{10}\\
w_7&=&(v_7+v_{11}+7v_{13}+v_{17}+92v_{19}-14v_{23}-994v_{29})/3^7\\
w_{11}&=&v_{11}/3^9-v_{13}/3^6-2v_{17}/3^8-13v_{19}/3^6+83v_{23}/3^8-11885v_{29}/3^9\\
w_{13}&=&v_{13}/3^5-v_{17}/3^6+7v_{19}/3^5+2v_{23}/3^6+227v_{29}/3^6\\
w_{17}&=&v_{17}/3^5+v_{19}/3^3+4v_{23}/3^5-26v_{29}/3^5\\
w_{19}&=&v_{19}/3^2-2v_{23}/3^3-v_{29}/3^3\\
w_{23}&=&(v_{23}-v_{29})/3^2\\
w_{29}&=&v_{29}.
\end{eqnarray*}
We chose to list here the basis in terms of eigenvectors $v_i$ rather than the classes $B_i$,
as was done in Propositions \ref{psiG2} to \ref{psiE8},
primarily for the sake of variety.

By performing a matrix multiplication,
we obtain the following formula for the Adams operations on the classes
$w_i$. This is totally analogous to Theorem \ref{psiE8}, except that here we are localized at 3.

\begin{prop}\label{e83ops} On the basis $\{w_i\}$ of $PK^1(E_8;\bz_{(3)})$ just defined, the Adams operations are given, for any integer $k$, by
\begin{eqnarray*}
\psi^k(w_1)&=&kw_1-\tfrac{730}3(k-k^7)w_7+(\tfrac{6560}3k-2190k^7+\tfrac{10}3k^{11})w_{11}\\&&
+(918k-\tfrac{24820}{27}k^7+\tfrac{10}9k^{11}+\tfrac4{27}k^{13})w_{13}\\
&&+(\tfrac{1484}3k-\tfrac{40150}{81}k^7+\tfrac{50}{81}k^{11}+\tfrac4{81}
k^{13}+\tfrac{28}{81}k^{17})w_{17}\\
&&+(\tfrac{121}3k-\tfrac{29200}{729}k^7
+\tfrac{10}{243}k^{11}-\tfrac{40}{729}k^{13}-\tfrac{28}{243}k^{17}
-\tfrac{109}{729}k^{19})w_{19}\\
&&+(-\tfrac{998}3k+\tfrac{727810}{2187}k^7-\tfrac{1030}{2187}k^{11}
-\tfrac{104}{2187}k^{13}-\tfrac{280}{2187}k^{17}-\tfrac{218}{2187}k^{19}
+\tfrac{1364}{2187}k^{23})w_{23}\\
&&+(\tfrac{8479}9k-\tfrac{6187480}{6561}k^7+\tfrac{99320}{59049}k^{11}
-\tfrac{316}{6561}k^{13}+\tfrac{364}{19683}k^{17}-\tfrac{109}{6561}k^{19}\\
&&+\tfrac{1364}{19683}k^{23}-\tfrac{44078}{59049}
k^{29})w_{29}\\
\psi^k(w_7)&=&k^7w_7-9(k^{7}-k^{11})w_{11}+(-\tfrac{34}9k^7+3k^{11}
+\tfrac79k^{13})w_{13}+(-\tfrac{55}{27}k^7+\tfrac53k^{11}\\
&&+\tfrac7{27}k^{13}
+\tfrac19k^{17})w_{17}
+(-\tfrac{40}{243}k^7+\tfrac19k^{11}-\tfrac{70}{243}k^{13}-\tfrac1{27}k^{17}+\tfrac{92}{243}k^{19})w_{19}\\
&&+(\tfrac{997}{729}k^7-\tfrac{103}{81}k^{11}-\tfrac{182}{729}k^{13}-\tfrac{10}{243}k^{17}+\tfrac{184}{729}k^{19}-\tfrac{14}{243}k^{23})w_{23}\\
&&+(-\tfrac{8476}{2187}k^7+\tfrac{9932}{2187}k^{11}-\tfrac{553}{2187}
k^{13}+\tfrac{13}{2187}k^{17}+\tfrac{92}{2187}k^{19}-\tfrac{14}{2187}k^{23}-\tfrac{994}{2187}k^{29})w_{29}\\
\psi^k(w_{11})&=&k^{11}w_{11}+\tfrac13(k^{11}-k^{13})w_{13}
+(\tfrac5{27}k^{11}-\tfrac19k^{13}-\tfrac2{27}k^{17})w_{17}\\
&&+(\tfrac1{81}k^{11}+\tfrac{10}{81}k^{13}+\tfrac2{81}k^{17}
-\tfrac{13}{81}k^{19})w_{19}\\
&&+(-\tfrac{103}{729}k^{11}+\tfrac{26}{243}
k^{13}+\tfrac{20}{729}k^{17}-\tfrac{26}{243}k^{19}+\tfrac{83}{729}k^{23})w_{23}\\
&&+(\tfrac{9932}{19683}k^{11}+\tfrac{79}{729}k^{13}-\tfrac{26}{6561}
k^{17}-\tfrac{13}{729}k^{19}+\tfrac{83}{6561}k^{23}-\tfrac{11885}
{19683}k^{29})w_{29}\\
\psi^k(w_{13})&=&k^{13}w_{13}+\tfrac13(k^{13}-k^{17})w_{17}
+(-\tfrac{10}{27}k^{13}+\tfrac19k^{17}+\tfrac7{27}k^{19})w_{19}\\
&&+(-\tfrac{26}{81}k^{13}+\tfrac{10}{81}k^{17}+\tfrac{14}{81}k^{19}
+\tfrac2{81}k^{23})w_{23}\\
&&+(-\tfrac{79}{243}k^{13}-\tfrac{13}{729}k^{17}
+\tfrac7{243}k^{19}+\tfrac2{729}k^{23}+\tfrac{227}{729}k^{29})w_{29}\\
\psi^k(w_{17})&=&k^{17}w_{17}-\tfrac13(k^{17}-k^{19 })w_{19}
+(-\tfrac{10}{27}k^{17}+\tfrac29k^{19}+\frac4{27}k^{23})w_{23}\\
&&+(\tfrac{13}{243}k^{17}+\tfrac1{27}k^{19}+\tfrac4{243}k^{23}
-\tfrac{26}{243}k^{29})w_{29}\\
\psi^k(w_{19})&=&k^{19}w_{19}+\tfrac23(k^{19}-k^{23})w_{23}
+(\tfrac19k^{19}-\tfrac2{27}k^{23}-\tfrac1{27}k^{29})w_{29}\\
\psi^k(w_{23})&=&k^{23}w_{23}+\tfrac19(k^{23}-k^{29})w_{29}\\
\psi^k(w_{29})&=&k^{29}w_{29}.
\end{eqnarray*}
\end{prop}

The second terms of the above expressions for $\psi^k(w_i)$ are, at the very least, strongly suggestive that $\Phi E_8$ can be built 
by fibrations from $\Phi S^n$'s according to the scheme in Diagram
\ref{att}. Here $\Phi$ is Bousfield's functor (\cite{Bous}) satisfying
$\pi_*(\Phi X)\approx\vp_*(X)$ (localized at 3). We do not wish to belabor this here, since it is not necessary for our analysis. However,
we note that it seems quite likely that this analysis could yield
information about the attaching maps between cells of the 
3-localizations of $E_8$ and $\Om E_8$, extending work in \cite{KM} and
\cite{HaHa}. Our results here were previously inaccessible, since
$\a_i$ is not detected by primary cohomology operations when $i>1$.

\begin{diag}\label{att}
\begin{center}
\begin{picture}(380,45)
\put(0,0){\framebox(20,20){3}}
\put(50,0){\framebox(20,20){15}}
\put(100,0){\framebox(20,20){23}}
\put(150,0){\framebox(20,20){27}}
\put(200,0){\framebox(20,20){35}}
\put(250,0){\framebox(20,20){39}}
\put(300,0){\framebox(20,20){47}}
\put(350,0){\framebox(20,20){59}}
\put(20,10){\line(1,0){30}}
\multiput(120,10)(50,0){5}{\line(1,0){30}}
\put(24,13){$3\a_3$}
\put(130,13){$\a_1$}
\put(180,13){$\a_2$}
\put(230,13){$\a_1$}
\put(280,13){$\a_2$}
\put(330,13){$\a_3$}
\put(60,20){\oval(90,10)[t]}
\put(55,28){$\a_5$}
\end{picture}
\end{center}
\end{diag}

We return now to the application of Proposition \ref{e83ops} to determining $\vp_*(E_8;3)$. By Theorem \ref{Bousthm},
$\vp_{2m}(E_8;3)$ has a presentation with 8 generators and 16 relations,
corresponding to $\psi^3(w_i)$ and $(\psi^2-2^m)(w_i)$, where we use
\ref{e83ops} for $\psi^2$ and $\psi^3$. By Proposition \ref{nu}, which we shall use frequently, the $\psi^2-2^m$ relations imply that 
$\vp_{2m}(E_8;3)=0$ if $m$ is even, and so we let $m=2k+1$.
We divide the $\psi^2-2^m$ relations by the unit 2, and arrange the
relations as rows of a $16\times 8$ matrix, with the $\psi^3$ relations listed first. 

We use {\tt Maple} to manipulate the matrix. 
The entry in position (15,8) is a unit, and so we pivot on it, 
and then delete the 15th row and the 
8th column. This corresponds to writing the 8th generator 
as a combination of the other generators, then removing that generator and 
the relation which expressed it in terms of the others, 
while substituting this relation into all the others.
In subsequent steps, we pivot on and then eliminate positions (14,7), (13,6), 
(12,5), (11,4), (1,3), and (8,2), ending up with a $9\times 1$ matrix
$G$. At each step, it is
essential that the element on which we pivot is a unit. It was by no means 
clear at the outset that this could be done 7 times, 
for it has the important consequence
that $\vp_{2m}(E_8;3)$ is cyclic for each integer $m$. 

In manipulating the matrix, we let, for $i\in\{0,6,10,12,16,18,22,28\}$, $P_i=2^{2k}-2^{i}$. Thus, for example, the elements in position (9,1) and (10,2) of the initial matrix will be $-P_0$ and $-P_6$.
After $j$ steps of pivoting ($j\le7$), all entries in the matrix will be quotients of polynomials in the $P_i$ of degree $\le j+1$. Since $P_i\equiv0$ mod 3, mod 3 values of polynomials are determined by their constant term,
which during the last few steps will be up to 30 digits long. Up to this point, we have only cared about mod 3 values, to find units on which we could pivot, but in the next step we need more delicate information about exponents of 3.

We will have $\vp_{4k+2}(E_8;3)\approx3^e$, where $e$ is the smallest exponent of 3 of the 
nine entries of the matrix $G$. The denominator polynomials are units, and can be ignored. 
{\tt Maple} does some factoring automatically. 
For example, the fourth entry of $G$ has a factor $3^{12}P_{10}P_{12}$, which we treat as $3^{14+\nu(k-5)+\nu(k-6)}$, since $\nu(P_{2j})=1+\nu(k-j)$. Throughout this section,
$\nu(-)=\nu_3(-)$.

Now we divide into cases depending upon the mod 9 value of $k$.
We consider first the case $k\equiv2\equiv11$ mod 9. After performing
the preliminary simplifications described in the preceding paragraph,
we replace each occurrence of $P_{2j}$ by $R+2^{22}-2^{2j}$. Thus
$R$ is representing $2^{2k}-2^{22}$, which satisfies $\nu(R)=1+\nu(k-11)$. We find that the nine relations are, up to unit multiples of all terms, as follows. We emphasize that this assumes that $k\equiv2$ mod 9,
so that, for example, $\nu(k-8)$ has been replaced by 1, and $\nu(k-10)$ 
by 0. We also point out that it seems to be infeasible to obtain these expressions by hand; their determination seems to require a computer. 
$$\begin{array}{l}
3^{25}+3^8R+3^6R^2+3^4R^3+3^2R^4+3R^5+R^6\\
3^2(3^{22}+3^7R+3^{8}R^2+3^3R^3+3^3R^4+3^2R^5+R^6)\\
3^{10}(3^{16}+3^6R+3^5R^2+3^2R^3+3R^4+R^5)\\
3^{15}(3^{12}+3^4R+3^2R^2+3R^3+R^4)\\
3^{22}(3^7+3^3R+3R^2+R^3)\\
3^{27}(3^3+3R+R^2)\\
3^{37+\nu(k-11)}\\
3^{24}+3^9R+3^7R^2+3^5R^3+3^3R^4+3^2R^5+3R^6+R^7\\
3^{10+\nu(k-11)}
\end{array}$$
Since $\nu(R)\ge3$, we find that the term $3^8R$ in the first relation
gives the smallest 3-power ($3^{9+\nu(k-11)}$) if $\nu(k-11)\le15$, while the term $3^{24}$ in the second or eighth relation will be smallest
if $\nu(k-11)\ge16$. This establishes the claim of Theorem \ref{e83thm}
for $\vp_{2m}(E_8;3)$ when $m=2k+1$ and $k\equiv 2$ mod 9.

The situation is similar when $k\equiv 8$ mod 9. If $R=2^{2k}-2^{16}$,
then the first two relations are, up to unit multiples of all terms,
$$3^{19}+3^8R+3^6R^2+3^4R^3+3^2R^4+3R^5+R^6$$
and
$$3^{18}+3^9R+3^8R^2+3^5R^3+3^5R^4+3^4R^5+3^2R^6,$$
while other relations involve larger exponents of 3. Note that these
refer to precisely the same relations as the first two of the nine listed above; however, $R$ now represents a different expression.
Since $\nu(R)\ge3$, the $R^i$-terms with $i\ge2$ are more highly
3-divisible than the $R^1$-term, and so may be ignored. We find that the smallest exponent of 3 is $\min(9+\nu(k-8),18)$, with the $9+\nu(k-8)$ coming from the first relation and the 18 from the second. This yields 
the case $k\equiv8$ mod 9 for $\vp_{2m}(-)$ in Theorem \ref{e83thm}.
The case $k\equiv 5$ mod 9 follows similarly by noting that if
$R=2^{2k}-2^{28}$, then the only significant terms are
$3^{32}+3^9R+3^6R^2$ in the first relation and $3^{30}+3^9R$ in the second and eighth. 
If $\nu(R)=3$, there was a possibility of cancellation of the $R^1$- and $R^2$-terms 
in the first relation, but not in the second.

The case $k\equiv 1$ mod 3 is easier. If $R=2^{2k}-2^8$, the first
relation is $3^6+3^5R+3^4R^2+3^3R^3+3^2R^4+3R^5+R^6$, which gives a relation $3^6$ if $\nu(R)\ge2$, while the other relations are more
highly 3-divisible. 

The case $k\equiv6$ mod 9 introduces a second-order effect. If $R=
2^{2k}-2^{12}$, the only significant terms are $3^{14}+3^6R$ from the
first relation and $3^{15}+3^9R$ from the second. 
The smallest 3-exponent is $7+\nu(k-6)$ if $\nu(k-6)<7$, and is
14 if $\nu(k-6)>7$. If $\nu(k-6)=7$, the two terms in the first relation are both $3^{14}$ 
times a unit, and we must analyze the mod 3 values of these units to tell whether the sum of 
these two terms is $3^{14}$ times a unit or is divisible by $3^{15}$. 
As the second relation is $3^{15}$ times a unit, we need only evaluate 
 the first relation mod $3^{15}$.

{\tt Maple} tells us that the unit coefficients of $3^{14}$ and $3^6R$ in the first relation
are both 1 mod 3. Thus mod $3^{15}$ the first relation is
$$3^{14}+3^6((1+3)^k-(1+3)^6)\equiv3^{14}+3^6\biggl(3(k-6)+3^2\bigl(\tbinom k2-\tbinom62\bigr)\biggr).$$
If $k=3^7u$, this becomes $3^{14}(1+u)$, and so is $3^{15}$ if $u\equiv
2$ mod 3, and $3^{14}$ if $u\equiv 1$. This yields the
$\min(7+\nu(k-6-2\cdot3^7),15)$ in \ref{e83thm} when $k\equiv 6$ mod 9.
The case $k\equiv0$ mod 9 follows similarly once {\tt Maple}
tells us that if $R=2^{2k}-2^{18}$, then the terms which can give
smallest 3-power are $(3a+1)3^{20}+(3b+2)3^6R$ from the first relation
and $3^{21}+3^9R$ from the second. 

The case $k\equiv 3$ mod 9 is easier. If $R=2^{2k}-2^6$, then the first relation begins
$3^8+3^7R+3^4R^2$, and other relations involve larger 3-powers. Thus $\vp_{2m}(-)\approx
\bz/3^8$ in this case, which appears in \ref{e83thm} bundled along with the case $k\equiv6$.

This concludes the computation of $\vp_{2m}(E_8;3)$ in Theorem \ref{e83thm}. The $\vp_{2m-1}(-)$ part of that theorem follows from
Propositions \ref{order} and \ref{number} and the following result.
\begin{prop}\label{onecyc} If $m=3$, then $\vp_{2m-1}(E_8;3)\approx\bz/3^6$.
\end{prop}
\begin{pf} By Proposition \ref{order} and the above computation that $\vp_6(E_8;3)\approx\bz/3^6$,
it suffices to show that $\vp_5(E_8;3)$ has an element $x$ satisfying $3^5x\ne0$.
For $k=2$ and 3, let $M_k$ denote the transpose of the matrix of $\psi^k$ with respect to the basis $\{w_1,w_7,w_{11},w_{13},
w_{17},w_{19},w_{23},w_{29}\}$ of Proposition \ref{e83ops}.
Thus $M_k$ is the matrix whose rows have as entries the various coefficients in the
equations of \ref{e83ops} with numerical values obtained by substituting
2 or 3 for $k$. Let $\Phi=M_2-8I$, an invertible matrix over $\bz_{(3)}$. Let $E=(1\ 1\ 1\ 1\ 1\ 1\ 1\ 1)$.
Then $EM_3$ is the sum of the rows of $M_3$, and this combination of the
$w_i$'s is in $\im(\psi^3)$.
Let $A=EM_3\Phi^{-1}$. {\tt Maple} computes 
\begin{eqnarray*}
A&=&(\tfrac{-1}2,\tfrac{548957}{120},\tfrac{-81625061}{2040},\tfrac
{-3845766033}{231880},\tfrac{-605269537661}{759778008},\tfrac{
1212671517798439}{488157370140},\\
&&\tfrac{4278408406308902221}{200144521757400},\tfrac{-248464589902573527757029703}{2238578548469645972700}).
\end{eqnarray*}
The second and third denominators are divisible by 3, but the others are not. 
Thus $3A$ represents an element of $PK^1(E_8;\bz_{(3)})$, and
$3A\Phi$, which represents $\psi^2-2^3$ applied to this vector,
is in $\im(\psi^3)$. 

Thus $3A\in\ker(\psi^2-2^3:PK^1(E_8;\bz_{(3)})/\im(\psi^3)\to PK^1(E_8;\bz_{(3)})/\im(\psi^3))$.
We will be done once we have shown that $3^6A\not\in\im(\psi^3)$, for then $3^5(3A)$ is a nonzero 
element of the kernel, which, by Theorem \ref{Bousthm}, is isomorphic to $\vp_5(E_8;3)$.
This is verified by {\tt Maple} by successively subtracting multiples of the rows of $M_3$ from
$3^6A$ to change the successive components of $3^6A$ to 0. For the first 7 steps, the required
multiplier will be in $\bz_{(3)}$, but when we get to the last step, our vector, which is equivalent
to $3^6A$ mod $\im(\psi^3)$, will be $(0,0,0,0,0,0,0,3^{28}u)$, with $u$ a unit in $\bz_{(3)}$.
Since the corresponding relation in $\im(\psi^3)$ is $3^{29}$, we conclude that $3^6A\not\in\im(\psi^3)$,
as desired.\end{pf}

This concludes the proof of Theorem \ref{e83thm}, with the only input being $\lambda^2$ in
$R(E_8)$. We find it useful to interpret our result in terms of homotopy charts such as
Diagram \ref{firstdiag}. We emphasize that this analysis is not part of our proof, but rather
is an attempt to understand how $\vp_*(E_8;3)$ is built from the $v_1$-periodic homotopy 
groups of the eight spheres  which build it. In \cite{E7}, in which the UNSS was the primary
tool, these charts were an integral part of the argument, but here we just use our computation
to see how those charts must have been filled in.  We are not saying that $E_8$ is built from
these eight spheres (It is not!), but rather that since $PK^1(E_8;\zphat)$ as an Adams module
is built from the related Adams modules for spheres (by our Theorem \ref{niceops}),
$\vp_*(E_8;3)$ is built from that of the spheres, and we could probably conclude that
Bousfield's $\Phi E_8$ is built from $\Phi$ applied to the spheres.

An important ingredient in the charts such as \ref{firstdiag} and those of \cite{E7} is the 
cyclicity of the groups for the \lq\lq sphere bundles'' described in Proposition \ref{sphbdl}.
In $(E_8,3)$, we also encounter \lq\lq sphere bundles'' with $\a_p$ as attaching map, a case
which was not considered in \ref{sphbdl}. Here we have the following result, whose 
proof we merely sketch.
\begin{prop}\label{ap} If $PK^1(X;\zphat)$ has generators $x$ and $y$ with $\psi^ky=k^{n+p(p-1)}
y$ and $\psi^kx=k^n+\frac u{p^\eps}k^n(k^{p(p-1)}-1)y$, where $u$ is a unit in $\bz_{(p)}$
and $\eps=1$ or $2$, then the only nonzero groups $\vp_*(X;p)$ are:
\footnote{If $\eps=2$ and $n=1$, then the coefficient of $y$ in $\psi^px$ is not
in $\bz_{(3)}$, and so this case is excluded.}

\noindent a. If $\eps=2$, then
$\vp_{2n+qi}(X;p)\approx\vp_{2n+qi-1}(X;p)\approx\bz/p^e$, where
$$e=\begin{cases}\min(n-1,2)&\text{if $\nu(i)=0$}\\
\min(n,\nu(i)+3)&\text{if $\nu(i-p)=1$}\\
\min(n+p(p-1),\nu(i-p)+3)&\text{if $\nu(i-p)>1$ and $n\ge4$}\\
\min(4+p(p-1),\nu(i-p-p^{p(p-1)+1})+n-1)&\text{if $\nu(i-p)>1$ and $2\le n\le3$}
\end{cases}$$

\noindent b. If $\eps=1$ and $n\ge2$, then $\vp_{2n+qi}(X;p)\approx
\vp_{2n+qi-1}(X;p)\approx\bz/p\oplus\bz/p^e$, where
$$e=\begin{cases}1&\text{if $\nu(i)=0$}\\
\min(n,\nu(i)+2)&\text{if $\nu(i-p)=1$}\\
\min(n+p(p-1),\nu(i-p)+2)&\text{if $\nu(i-p)>1$ and $n\ge3$}\\
\min(3+p(p-1),\nu(i-p-p^{p(p-1)+1})+1)&\text{if $\nu(i-p)>1$ and $n=2$}
\end{cases}$$

\noindent c. If $\eps=1$ and $n=1$, then 
$$\vp_{2n+qi}(X;p)\approx\vp_{2n+qi-1}(X;p)
\approx\bz/p^{\min(2+p(p-1),\nu(i-p)+1)}.$$
\end{prop}
\begin{pf} The proof follows closely that of Proposition \ref{sphbdl}
with $t=p$. We let $s=r^{p-1}=1+\b p$, and note $\nu(s^i-1)=1+\nu(i)$.
We consider first the case $\vp_{2n+qi}(X;p)$.

First, let $\eps=2$. As occurred in \ref{sphbdl}, the fourth relation
allows one to eliminate $y$ and deduce that the group is cyclic with relations (on generator $x$)
$$p^{1+\nu(i)}p^{\min(n+p(p-1),\nu(i-p)+1)}\text{\quad and\quad}p^{n-2}(p^{p(p-1)}(s^i-1)
+s^p-s^i).$$
The only case in which it is not easily seen that the minimal exponent is as claimed occurs 
when $n<4$ and $i-p=u_0p^{p(p-1)+1}$
with $u_0$ a unit in $\bz_{(p)}$. In this case, the first relation is
$p^{p(p-1)+4}$, and so we analyze the second relation mod $p^{p(p-1)+4}$,
obtaining
\begin{eqnarray*}
&&p^{n-2}\biggl(p^{p(p-1)}\bigl((1+\b p)^{i}-1\bigr)-(1+\b p)^p\bigl((1+\b p)^{i-p}-1\bigr)\biggr)\\
&\equiv&\b p^{n+p(p-1)}(1-u_0)(1-\b p/2),
\end{eqnarray*}
which yields the claim of the proposition.

Now let $\eps=1$ and $n=1$. We use the second of the four relations similar to those in the proof of 
\ref{sphbdl} to eliminate
$y$, again obtaining that $\vp_{2n+qi}(X;p)$ is cyclic with relations on
generator $x$
$$p^{p(p-1)+2},\ p(s^{i-p}-1), \text{ and
}p^{p(p-1)}(s^i-1)-s^p(s^{i-p}-1),$$
whose minimal exponent of $p$ is easily seen to be $\min(2+p(p-1),1+\nu(i-p))$, as claimed.

Finally, let $\eps=1$ and $n>1$. All terms in all four relations are divisible by $p$. We split off a
$\bz/p$ generated by $1/p$ times the last relation, and replace $uy$ by $(p(s^i-1)/(s^p-1))x$.
The remaining summand has relations (on $x$)
$$p^{n+p(p-1)+\nu(i)},\ p^{\nu(i)+\nu(i-p)+1},\text{ and }p^{n-2}(p^{p(p-1)}(s^i-1)+s^p-s^i).$$
The only case in which it is not easily seen that the smallest exponent is as claimed occurs 
when $i-p=u_0p^{p(p-1)+1}$ with $u_0$ a unit in $\bz_{(p)}$. Similarly to the case $\eps=2$,
the last relation becomes $\b p^{2+p(p-1)}(1-u_0)$ mod $p^{p(p-1)+3}$,
which yields the claim of the 
proposition.

The result for $\vp_{2n+qi-1}(X;p)$ when $\eps=2$ follows from the
calculation of $v_{2n+qi}(X;p)$
just completed together with Propositions \ref{order}, \ref{number}, and \ref{case}.
\begin{prop}\label{case} If $\eps=2$ and $\nu(i-p)\ge p(p-1)+n$, then $\vp_{2n+qi-1}(X;p)$
is cyclic.\end{prop}
\begin{pf} The proof is similar to that of Proposition \ref{sphbdl}. In the case considered here,
$|\vp_{2n+qi-1}(X)|=p^{p(p-1)+n}$ and, by (\ref{sph}), $\vp_{2n+qi-1}(S^{2n+pq+1})\approx
\bz/p^{p(p-1)+n}$, and so the result follows from the exact sequence
$$\to\vp_{2n+qi-1}(S^{2n+1})\to\vp_{2n+qi-1}(X)\to\vp_{2n+qi-1}(S^{2n+pq+1})\to 0.$$
\end{pf}

The result for $\vp_{2n+qi-1}(X)$ when $\eps=1$ and $n=1$ follows from Propositions
\ref{order} and \ref{number} and the fact that if $i\not\equiv0$ mod $p$, then $\vp_{2n+qi}(X)$
is cyclic, since it has order $p$.

The proof for $\vp_{2n+qi-1}(X)$ when $\eps=1$ and $n\ge2$ is more difficult, since it is not
enough to verify it for one value of $i$. We focus on the case $\nu(i-p)=1$, $\nu(i)\le n-2$.
Other cases are handled similarly. By Theorem \ref{Bousthm}
$$\vp_{2n+qi-1}(X)\approx\ker(\phi:G\to G),$$
where
\begin{eqnarray*}
G&=&\langle x,y:\ p^{n+p(p-1)}y,\ p^nx+up^{n-1}(p^{p(p-1)}-1)y\rangle\\
\phi(y)&=&r^n(s^p-s^{i})y\\
\phi(x)&=&r^n((1-s^{i})x+\tfrac up(s^p-1)y).
\end{eqnarray*}
Using the second relation in $G$, we compute
$$\phi(p^{n-\nu(i)-1}x)=r^nup^{n-\nu(i)-2}(s^p-s^{i})(1+Ap^{p(p-1)})y,$$
where $A:=(s^{i}-1)/(s^p-s^{i})$ satisfies $\nu(A)=\nu(i)-1$. Then
$$z:=p^{n-\nu(i)-1}x-up^{n-\nu(i)-2}(1+Ap^{p(p-1)})y$$
is in $\ker(\phi)$, and $p^{\nu(i)+1}z=-up^{n-1+p(p-1)}u'y$, 
where $u'=(s^p-1)/(s^p-s^{i})$ is a unit. Thus this element $p^{\nu(i)+1}z$ has order $p$ in $G$,
and so $\ker(\phi)\approx\bz/p^{\nu(i)+2}\oplus\bz/p$, with generators $z$ and 
$uu'p^{n+p(p-1)-2}y+p^{\nu(i)}z$.
\end{pf}

We list in Diagrams \ref{ch1} and \ref{ch2} some representative charts indicating how
$\vp_{4k+\delta}(E_8)$, $\delta=1,2$, is built from the various $\vp_{4k+\delta}(\sn)$.
These charts are of the type of Diagram \ref{firstdiag} and those of \cite{E7}.
For each of the four cases considered in each diagram, the left tower is building
$\vp_{4k+1}(E_8)$ and the right tower is building $\vp_{4k+2}(E_8)$. The numbers on the 
left side of each diagram are the dimensions of the spheres. Dots represent $\bz/3$, and
a number $e$ represents $\bz/p^e$. Vertical lines are extensions (multiplication by 3), as is
the curved line from 23 to 3. Slanting arrows are boundary morphisms in exact sequences.
In the third and fourth case of both diagrams, the boundary from the bottom generator
is hitting a sum of two classes, which is relevant to the claim that the cokernel is cyclic. 

The extensions from 59 to 47, and from 15 to 3, are consequences of Proposition \ref{ap}
together with the terms $\frac{730}3(k-k^7)w_7$ and $\frac19(k^{23}-k^{29})w_{29}$ in Proposition
\ref{e83ops}. The other extensions are derived similarly from Proposition \ref{sphbdl}.
 The most unexpected part of the charts is the extension in Diagram \ref{ch2}
from 27 into the sum of 23 and 15 in dimension $4k+2$. 
This was seen to be necessary in order that $\vp_{4k+2}(E_8)$
be cyclic in these cases. A check that this extension is really present was made by having
{\tt Maple} compute $v_{4k+2}(X(15,23,27))$, corresponding to the indicated subquotient Adams 
module from \ref{e83ops}. This was done with $k=30$, corresponding to the first case in 
\ref{ch2}, and $\bz/3\oplus\bz/3^6$ was obtained, consistent with the unexpected extension.
This extension is probably also present in the cases of Diagram \ref{ch1}, but in those cases it 
would not affect the result. 

There is also an unexpected extension in $4k+1$ from 23 to $p$ times the generator in 15
when $k\equiv 6$ mod 9. This seemed necessary in order to get the correct answer when
$\nu(k-6)\ge12$, and was confirmed by a {\tt Maple} computation that $\vp_{4k+2}(X(15,23,27))
\approx\bz/3\oplus\bz/3^{\min(3+\nu(k-6),13)}$ if $k\equiv6$ mod 9. The orders of the 
cyclic groups in the eight cases below are 6, $\nu(k-11)+9$, 24, 24, 8, 8, $\nu(k-6)+7$, and
14, in agreement with Theorem \ref{e83thm}.

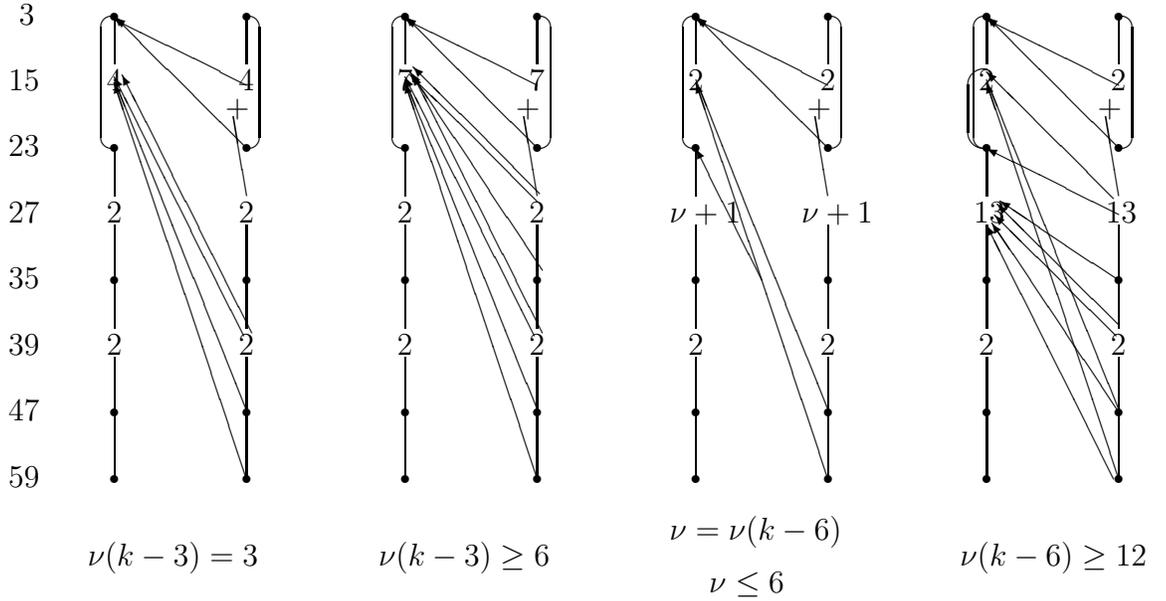
\begin{diag}\label{ch1}

\begin{center}
\begin{picture}(430,230)
\def\elt{\circle*{3}}
\def\mp{\multiput}
\put(0,40){$59$}
\put(0,65){$47$}
\put(0,90){$39$}
\put(0,115){$35$}
\put(0,140){$27$}
\put(0,165){$23$}
\put(0,190){$15$}
\put(4,215){$3$}
\put(40,10){$k\equiv1\ (3)$}
\put(148,0){$2\le\nu\le15$}
\put(142,20){$\nu=\nu(k-11)$}
\put(245,10){$\nu(k-11)=16$}
\put(355,10){$\nu(k-11)\ge22$}
\mp(40,43)(0,25){8}{\elt}
\mp(90,43)(0,25){8}{\elt}
\mp(40,43)(50,0){2}{\line(0,1){125}}
\mp(40,193)(110,0){4}{\line(0,1){25}}
\mp(90,193)(110,0){4}{\line(0,1){25}}
\mp(40,193)(110,0){4}{\oval(10,50)[l]}
\mp(90,193)(110,0){4}{\oval(10,50)[r]}
\mp(147,40)(110,0){3}{$2$}
\mp(197,40)(110,0){3}{$2$}
\mp(147,115)(110,0){3}{$2$}
\mp(197,115)(110,0){3}{$2$}
\mp(147,165)(110,0){3}{$2$}
\mp(197,165)(110,0){3}{$2$}
\mp(150,93)(110,0){3}{\elt}
\mp(200,93)(110,0){3}{\elt}
\mp(150,143)(110,0){3}{\elt}
\mp(200,143)(110,0){3}{\elt}
\mp(150,193)(110,0){3}{\elt}
\mp(200,193)(110,0){3}{\elt}
\mp(150,218)(110,0){3}{\elt}
\mp(200,218)(110,0){3}{\elt}
\mp(140,65)(50,0){2}{$\nu+1$}
\mp(256,65)(50,0){2}{$17$}
\mp(364,65)(50,0){2}{$23$}
\mp(150,49)(110,0){3}{\line(0,1){15}}
\mp(200,49)(110,0){3}{\line(0,1){15}}
\mp(150,74)(110,0){3}{\line(0,1){40}}
\mp(200,74)(110,0){3}{\line(0,1){40}}
\mp(150,124)(110,0){3}{\line(0,1){40}}
\mp(200,124)(110,0){3}{\line(0,1){40}}
\mp(90,43)(110,0){4}{\vector(-1,3){50}}
\mp(90,193)(110,0){4}{\vector(-2,1){50}}
\mp(306,40)(1,4){2}{\vector(-1,3){42}}
\mp(418,44)(0,3){2}{\vector(-2,1){46}}
\put(418,67){\vector(-2,1){48}}
\mp(418,67)(0,3){2}{\vector(-1,1){46}}
\put(418,70){\vector(-2,3){48}}
\mp(418,70)(1,3){2}{\vector(-1,2){46}}
\end{picture}
\end{center}
\end{diag}

\begin{diag}\label{ch2}
\begin{center}
\begin{picture}(425,225)
\def\elt{\circle*{3}}
\def\mp{\multiput}
\put(0,40){$59$}
\put(0,65){$47$}
\put(0,90){$39$}
\put(0,115){$35$}
\put(0,140){$27$}
\put(0,165){$23$}
\put(0,190){$15$}
\put(4,215){$3$}
\put(30,10){$\nu(k-3)=3$}
\put(140,10){$\nu(k-3)\ge6$}
\put(250,20){$\nu=\nu(k-6)$}
\put(265,0){$\nu\le6$}
\put(360,10){$\nu(k-6)\ge12$}
\mp(40,43)(110,0){4}{\elt}
\mp(90,43)(110,0){4}{\elt}
\mp(40,68)(110,0){4}{\elt}
\mp(90,68)(110,0){4}{\elt}
\mp(40,118)(110,0){4}{\elt}
\mp(90,118)(110,0){4}{\elt}
\mp(40,168)(110,0){4}{\elt}
\mp(90,168)(110,0){4}{\elt}
\mp(40,218)(110,0){4}{\elt}
\mp(90,218)(110,0){4}{\elt}
\mp(37,90)(110,0){4}{$2$}
\mp(87,90)(110,0){4}{$2$}
\mp(37,140)(110,0){2}{$2$}
\mp(87,140)(110,0){2}{$2$}
\mp(37,190)(50,0){2}{$4$}
\mp(147,190)(50,0){2}{$7$}
\mp(257,190)(50,0){2}{$2$}
\mp(367,190)(50,0){2}{$2$}
\mp(250,140)(50,0){2}{$\nu+1$}
\mp(365,140)(50,0){2}{$13$}
\mp(40,43)(110,0){4}{\line(0,1){46}}
\mp(90,43)(110,0){4}{\line(0,1){46}}
\mp(40,100)(110,0){4}{\line(0,1){39}}
\mp(90,100)(110,0){4}{\line(0,1){39}}
\mp(40,200)(110,0){4}{\line(0,1){18}}
\mp(90,200)(110,0){4}{\line(0,1){18}}
\mp(40,193)(110,0){4}{\oval(10,50)[l]}
\mp(90,193)(110,0){4}{\oval(10,50)[r]}
\mp(40,150)(110,0){4}{\line(0,1){18}}
\mp(82,180)(110,0){4}{$+$}
\mp(90,150)(110,0){4}{\line(-1,6){5}}
\mp(90,192)(110,0){4}{\vector(-2,1){50}}
\mp(90,168)(110,0){4}{\vector(-1,1){50}}
\mp(90,43)(110,0){4}{\vector(-1,3){50}}
\mp(90,69)(110,0){4}{\line(-2,5){50}}
\mp(89,97)(110,0){2}{\vector(-1,2){49}}
\mp(92,98)(110,0){2}{\vector(-1,2){49}}
\put(202,122){\vector(-2,3){50}}
\mp(200,148)(1,3){2}{\vector(-1,1){48}}
\put(285,118){\vector(-1,2){25}}
\put(370,183){\oval(14,30)[l]}
\put(418,149){\vector(-1,1){48}}
\put(420,143){\vector(-2,1){50}}
\put(420,118){\vector(-3,2){45}}
\mp(419,97)(1,4){2}{\vector(-1,1){46}}
\put(420,68){\vector(-2,3){48}}
\put(418,43){\vector(-1,2){48}}
\end{picture}
\end{center}
\end{diag}

The charts when $k\equiv5$,8 mod 9 are very similar to those when $k\equiv 2$ given in Diagram
\ref{ch1}, while those when $k\equiv0$ mod 9 are very similar to those when $k\equiv 6$ given in
Diagram \ref{ch2}. Again we emphasize that these charts are not a part of our proof, but rather
a way of interpreting our result in a way which is closer to past methods of calculating
homotopy groups.

\section{{\tt LiE} program for computing $\lambda^2$ in $R(E_8)$}\label{pgm}
In this section we describe the program written in the specialized software
{\tt LiE} (\cite{Lie})
to perform the calculation described in Section \ref{rep}.
We list the program and then describe what it is doing.

{\obeylines\tt
setdefault E8
on + height
mm=id(8)
for r row mm do
\quad ext=alt\_tensor(2,r); extt=ext;
\quad p2=0X null(8); x=1;
\quad while x==1 do x=0;
\qq  for i=1 to length(extt) do u=expon(extt,i);
\qq\quad   if u[1]+u[2]+u[3]+u[4]+u[5]+u[6]+u[7]+u[8]>1 then j=1;
\qq\qq    while u[j]==0 do j=j+1 od;
\qq\qq    v=null(8); v[j]=u[j]; w=u-v;
\qq\qq if w==null(8) then w[j]=1; v[j]=v[j]-1 fi;
\qq\qq    p1=tensor(v,w); n=length(p1); utop=expon(p1,n);
\qq\qq    if extt | v==0 then x=1; p2=p2+1X v fi;
\qq\qq    if extt | w==0 then x=1; p2=p2+1X w fi;
\qq\qq    for k=1 to n-1 do a=expon(p1,k);
\qq\qq\quad     if extt | a==0 then x=1; p2=p2+1X a fi od
\qq\quad    fi
\qq   od;
\qq   extt=extt+p2
\quad  od;
\quad pdim=0X null(8); pder1=1X[1,0,0,0,0,0,0,0];
\quad pder2=1X[0,1,0,0,0,0,0,0]; pder3=1X[0,0,1,0,0,0,0,0];
\quad pder4=1X[0,0,0,1,0,0,0,0]; pder5=1X[0,0,0,0,1,0,0,0];
\quad pder6=1X[0,0,0,0,0,1,0,0]; pder7=1X[0,0,0,0,0,0,1,0];
\quad pder8=1X[0,0,0,0,0,0,0,1];
\quad for i=1 to length(extt) do
\qq u=expon(extt,i); pdim=pdim+dim(u)X u;
\qq  if u[1]+u[2]+u[3]+u[4]+u[5]+u[6]+u[7]+u[8]>1 then j=1;
\qq\quad   while u[j]==0 do j=j+1 od;
\qq\quad   v=null(8); v[j]=u[j]; w=u-v;
\qq\quad if w==null(8) then w[j]=1; v[j]=v[j]-1 fi;
\qq\quad   p1=tensor(v,w); n=length(p1);
\qq\quad   c1=dim(v)*(pder1 | w)+dim(w)*(pder1 | v);
\qq\quad   c2=dim(v)*(pder2 | w)+dim(w)*(pder2 | v);
\qq\quad   c3=dim(v)*(pder3 | w)+dim(w)*(pder3 | v);
\qq\quad    c4=dim(v)*(pder4 | w)+dim(w)*(pder4 | v);
\qq\quad   c5=dim(v)*(pder5 | w)+dim(w)*(pder5 | v);
\qq\quad   c6=dim(v)*(pder6 | w)+dim(w)*(pder6 | v);
\qq\quad   c7=dim(v)*(pder7 | w)+dim(w)*(pder7 | v);
\qq\quad   c8=dim(v)*(pder8 | w)+dim(w)*(pder8 | v);
\qq\quad   for i=1 to n-1 do c1=c1-coef(p1,i)*(pder1 | expon(p1,i));
\qq\qq     c2=c2-coef(p1,i)*(pder2 | expon(p1,i));
\qq\qq     c3=c3-coef(p1,i)*(pder3 | expon(p1,i));
\qq\qq     c4=c4-coef(p1,i)*(pder4 | expon(p1,i));
\qq\qq     c5=c5-coef(p1,i)*(pder5 | expon(p1,i));
\qq\qq     c6=c6-coef(p1,i)*(pder6 | expon(p1,i));
\qq\qq     c7=c7-coef(p1,i)*(pder7 | expon(p1,i));
\qq\qq     c8=c8-coef(p1,i)*(pder8 | expon(p1,i)) od;
\qq\quad   pder1=pder1+c1 X u; pder2=pder2+c2 X u;
\qq\quad   pder3=pder3+c3 X u; pder4=pder4+c4 X u;
\qq\quad   pder5=pder5+c5 X u; pder6=pder6+c6 X u;
\qq\quad pder7=pder7+c7 X u;   pder8=pder8+c8 X u
\qq  fi
\quad od;
\quad der1=0; der2=0; der3=0; der4=0; der5=0; der6=0; der7=0; der8=0;
\quad for i=1 to length(ext) do
\qq der1=der1+coef(ext,i)*(pder1 | expon(ext,i));
\qq  der2=der2+coef(ext,i)*(pder2 | expon(ext,i));
\qq  der3=der3+coef(ext,i)*(pder3 | expon(ext,i));
\qq  der4=der4+coef(ext,i)*(pder4 | expon(ext,i));
\qq  der5=der5+coef(ext,i)*(pder5 | expon(ext,i));
\qq  der6=der6+coef(ext,i)*(pder6 | expon(ext,i));
\qq  der7=der7+coef(ext,i)*(pder7 | expon(ext,i));
\qq  der8=der8+coef(ext,i)*(pder8 | expon(ext,i))
\quad od;
\quad print(der1);print(der2);print(der3);print(der4);
\quad print(der5);print(der6);print(der7);print(der8)
od}

The first line says that the program is always working with the Lie algebra
$E_8$. The second line says that monomials are ordered by increasing height,
a notion which was defined in Section \ref{rep}. The first {\tt for} loop,
which extends throughout the program, lets {\tt r} run over rows of an
$8\times 8$ identity
matrix, with the $i$th iteration
corresponding to computing $\lambda^2(\rhot_i)$. The variable {\tt ext} is
$\lambda^2(\rho_i)$ written as a polynomial whose exponents are
dominant weights, and coefficients are their multiplicity. The variable
{\tt extt} will be a modified version of {\tt ext},
expanded to include additional
terms whose derivatives must be computed, as described in Section \ref{rep}.
The variable {\tt x} tells whether any new terms were adjoined to {\tt extt}
in the most recent iteration. Each exponent {\tt u} in {\tt extt} with sum of
entries greater than 1 is decomposed as {\tt v}$+${\tt w}.  The polynomial
{\tt p1} represents
$V({\bold v})\ot V(\bold w)$. If {\tt v} or {\tt w} or any term
of {\tt p1} does not appear in {\tt extt}, then it is adjoined to {\tt extt},
because we will need to know its derivatives. The purpose of the 14-line
{\tt while} loop is just to adjoin these terms.

The variable {\tt pdim} is a polynomial whose coefficient of {\tt X u}
is $\dim(V(\bold u))$, while variables {\tt pder}$j$ are polynomials whose
coefficient of {\tt X u} is $\partial_j(V(\bold u))$,
the coefficient of $\rhot_j$ in $L(V(\bold u))$. Here {\tt u} ranges over
all terms in {\tt extt}; these will depend upon which $\lambda^2(\rhot_i)$ we
are computing. The long {\tt for} loop computes these inductively. It still
works with {\tt extt}, which is $\lambda^2(\rhot_i)$ modified to include extra
terms whose derivatives are needed in the induction. The variable {\tt c}$j$
is $\partial_j(V(\bold u))$, i.e., the coefficient of $\rhot_j$ in $V(\bold u)$.
It is determined by writing ${\bold u}=\bold v+\bold w$, expanding
$$V({\bold v})\ot V({\bold w})=V({\bold u})+\sum k_lV({\bold t}_l)$$
with height(${\bold t}_l)<$height($\bold u$),
so that $\partial_j(V({\bold t}_l))$ has
already been computed, and computing
$$\partial_j(V({\bold u}))=\dim(V({\bold v}))\partial_j(V({\bold w}))+
\dim(V({\bold w}))\partial_j(V({\bold v}))
-\sum k_l\partial_j(V({\bold t}_l)).$$

Finally, in the last short {\tt for} loop, $\partial_j(\lambda^2(\rhot_i))$,
represented by {\tt der}$j$, is computed by forming the appropriate combination
of the derivatives of the terms of $\lambda^2(\rhot_i)$. These terms are in the
polynomial {\tt ext}, and the derivatives are
incorporated as the coefficients in {\tt pder}$j$.

This program does not perform the subtraction of $\dim(\rho_i)\rhot_i$, which
is the last step in obtaining $L(\lambda^2\rhot_i)$. (See (\ref{14}).)
It was more convenient to do this in the {\tt Maple} program which was the next
step in our analysis.
No doubt, this {\tt LiE} program could be written more efficiently with arrays,
but the author felt more comfortable with it in the given form.

\section {Analysis of $F_4$ and $E_7$ at the prime 3}\label{F4}
When the analysis of Sections \ref{rep} and \ref{vecs} is performed for $F_4$
at the prime 3, we find that its Adams operations are isomorphic to those of
$B(11,15)\times B_5(3,23)$, where the second factor is an $S^3$-bundle over
$S^{23}$ with attaching map $\a_5$. Since $F_4$ is known to be 3-equivalent
to $B(11,15)\times H$, where $H$ is Harper's torsion $H$-space, this suggests
that there is an isomorphism of Adams modules
$$PK^1(H)_{(3)}\approx PK^1(B_5(3,23))_{(3)}$$
at the prime 3, or, more generally,
$$PK^1(H_p)_{(p)}\approx PK^1(B_{p+2}(3,3+q(p+2)))_{(p)},$$
for any odd prime $p$. This is implied by a result of Kono (\cite{Kono})
on Chern character in $K^*(H_p)$. However, its implication for
$v_1$-periodic homotopy groups does not quite agree with those determined
for $H_3$ in \cite{F4} or $H_p$ in \cite{Tel}.

The groups $\vp_*(B_{p+2}(3,3+q(p+2)))$ computed in Proposition \ref{sphbdl}
agree with those computed by the UNSS in \cite[2.1]{BDMi}. They agree for the
most part with the groups $\vp_*(H_p)$ determined in \cite[3.2]{Tel}. The
calculation in Proposition \ref{sphbdl} gives
$$\vp_{qi+2}(H_p)\approx\vp_{qi+1}(H_p)\approx\bz/p^{\min(p^2+p,1+\nu_p(i-p-2-\a
p^{p^2+p-2}))},$$
with $\a=p+2$,
while that in \cite[3.2]{Tel} is the same with $\a=1$. This discrepancy
led the author to find
a small mistake in the delicate argument presented in \cite{Tel} (for any odd
prime $p$) and in \cite{F4} (for $p=3$). The mistake in \cite{F4} occurs in
the fifth bulleted item on page 300. The correct exponent of $v_1$ there is
not 3, but rather a large number $e$ which is not a multiple of 3. Then
$(\eta_R(v_1^e)-v_1^e)\ot h_2^3$ does not desuspend to $S^5$, as would have
been the case if $e=3$. This affects \cite[(2.18)]{F4}, causing the 1 at the
end of it to be changed to a 2. The same mistake occurs in \cite{Tel} in the
second half of page 91; in this case the exponent of $v_1$ was overlooked
entirely. This was one of the rare places where exponents of $v_1$ are
consequential.

We have used the methods of this paper to check the result for $\vp_{\text{ev}}
(E_7;3)$ obtained in \cite{E7}. We obtain complete agreement with the results
presented there. In that paper, it was stated that if $j$ is odd and
$j\equiv5$ or 8 mod 9, then
$$\vp_{2j}(E_7)\approx
\vp_{2j-1}(E_7)\approx\bz/p\oplus\bz/3^{\min(19,\nu(j-17-2\delta\cdot3^{13})+4)
},$$
with $\delta=2$, 5, or 8. The calculation based upon the methods of this paper
shows that $\delta=5$. Also, in \cite{E7} it was stated that if $j$ is odd and
$\nu_3(j-11)\ge10$, then
$$\vp_{2j}(E_7)\approx\bz/3^3\oplus\bz/3^{12}\text{ or
}\bz/3^4\oplus\bz/3^{11}.$$
By the methods of this paper, we can show that the first of these splittings
is the valid one.

\def\line{\rule{.6in}{.6pt}}

\end{document}